\title{Reduced basis approximation and \emph{a posteriori} error estimation: applications to elasticity problems in several parametric settings}
\author{
        Dinh Bao Phuong Huynh \\
         Akselos SA and Department of Mechanical Engineering \\
         Massachusetts Institute of Technology, room 3-264 \\
          77 Mass Avenue, MA02142, Cambridge, USA\\
         \and
        Federico Pichi\\
        mathLab, Mathematics Area, SISSA,\\ International School for Advanced Studies \\           	    via Bonomea 265, 34136 Trieste, Italy\\
        \and
        Gianluigi Rozza \\
        mathLab, Mathematics Area, SISSA, \\International School for Advanced Studies \\           	    via Bonomea 265, 34136 Trieste, Italy\\
}
\date{\today}

\documentclass[12pt]{article}
\usepackage{mathptmx}       
\usepackage{helvet}         
\usepackage{courier}        
\usepackage{type1cm}        
\usepackage{makeidx}         
\usepackage{graphicx}        
\usepackage{multicol}        
\usepackage[bottom]{footmisc}

\usepackage{amssymb}
\usepackage{amsmath}
\usepackage{graphicx}
\usepackage{epsfig}
\usepackage{subfigure}
\usepackage{multirow}
\usepackage{color}
\usepackage{newtxtext,newtxmath}

\graphicspath{{Figures/}}

\newcommand{\refeq}[1]{(\ref{#1})}

\newcommand{\bfzeta}{\boldsymbol{\zeta}}

\newcommand{\calN}{\mathcal{N}}
\newcommand{\calE}{\mathcal{E}}
\newcommand{\p}{\partial}
\def\bfmu{\boldsymbol{\mu}}
\def\bfchi{\boldsymbol{\chi}}
\def\bfpsi{\boldsymbol{\psi}}
\def\bfx{\mathbf{x}}
\def\bfu{\mathbf{u}}

\def\bfw{\mathbf{w}}
\def\bfr{\mathbf{r}}
\def\bfy{\mathbf{y}}
\def\bfehat{\mathbf{\hat{e}}}

\def\bfP{\mbox{\boldmath$P$}}
\def\bfT{\mbox{\boldmath$T$}}

\def\bfB{\mathbf{B}}
\def\bfG{\mathbf{G}}
\def\bfP{\mathbf{P}}

\def\bfD{\mathbf{D}}
\def\bfH{\mathbf{H}}
\def\bfR{\mathbf{R}}
\def\bfS{\mathbf{S}}
\def\bfE{\mathbf{E}}
\def\bfK{\mathbf{K}}
\def\bfF{\mathbf{F}}
\def\bfL{\mathbf{L}}
\def\bfY{\mathbf{Y}}
\def\bfZ{\mathbf{Z}}
\def\bfC{\mathbf{C}}
\def\bfe{\mathbf{e}}
\def\calD{\mbox{\boldmath$\mathcal{D}$}}

\def\bfehat{\mbox{\boldmath$\hat{e}$}}

\begin{document}
\maketitle

\begin{abstract}
In this work we consider (hierarchical, Lagrange) reduced basis approximation and {\em a posteriori\/} error estimation for elasticity problems in affinley parametrized  geometries.  The essential ingredients of the methodology are: a Galerkin projection onto a low-dimensional space associated with a smooth ``parametric
manifold'' --- dimension reduction; an efficient and effective
greedy sampling methods for identification of optimal and
numerically stable approximations --- rapid convergence;
an {\em a posteriori\/} error estimation procedures --- rigorous and
sharp bounds for the functional outputs related with the underlying solution or related quantities of interest, like stress intensity factor; and 
Offline-Online computational decomposition strategies --- minimum
{\em marginal cost\/} for high performance in the real-time and many-query (e.g., design and optimization) contexts. We
present several  illustrative results for  linear elasticity problem in parametrized geometries representing 2D Cartesian or 3D axisymmetric configurations like an arc-cantilever beam, a center crack problem, a composite  unit cell or a woven composite beam, a multi-material plate, and a closed vessel.
We consider different parametrization for the systems: either physical quantities --to model the materials and loads-- and geometrical
parameters --to model different geometrical configurations-- with isotropic and orthotropic
materials working in plane stress and plane strain approximation. We would like to underline the versatility of the methodology in very different problems. As last example we provide a nonlinear setting with increased complexity.
\end{abstract}

\section{Introduction}
\label{sec:1}
In several fields, from continuum mechanics to fluid dynamics, we need to solve numerically very complex problems that arise from physics laws. Usually we model these phenomena through partial differential equations (PDEs) and we are interested in finding the field solution and also some other quantities that increase our knowledge on the system we are describing. Almost always we are not able to obtain an analytical solution, so we rely on some discretization techniques, like Finite Element (FE) or Finite Volume (FV), that furnish an approximation of the solution. We refer to this methods as the ``truth'' ones, because they require very high computational costs, especially in parametrized context.
In fact if the problem depends on some physical or geometrical parameter, the \emph{full-order} or \emph{high-fidelity} model has to be solved many times and this might be quite demanding. Examples of typical applications of relevance are optimization, control, design, bifurcation detection and real time query.
For this class of problems, we aim to replace the high-fidelity problem by one of much lower numerical complexity, through the \emph{model order reduction} approach \cite{MOR2016}.
We focus on Reduced Basis (RB) method \cite{hesthaven2015certified,Quarteroni2011,quarteroni2015reduced,morepas2017,benner2017model} which provides both \emph{fast} and \emph{reliable} evaluation of an input (parameter)-output relationship. The main features of this methodology are \emph{(i)} those related to the classic Galerkin projection on which RB method is built upon \emph{(ii)} an \emph{a posteriori} error estimation which provides sharp and rigorous bounds and \emph{(iii)} offline/online computational strategy which allows rapid computation.
The goal of this chapter is to present a very efficient \emph{a posteriori} error estimation for linear elasticity parametrized problem.  We show many different configurations and settings, by applying RB method to approximate problems using plane stress and plane strain formulation and to deal both with isotropic and orthotropic materials. We underline that the setting for very different problems is the same and unique.

This work is organized as follows. In Section 2, we first present a ``unified'' linear elasticity formulation; we then briefly introduce the geometric mapping strategy based on domain decomposition; we end the Section with the affine decomposition forms and the definition of the ``truth'' approximation, which we shall build our RB approximation upon. In Section 3, we present the RB methodology and the offline-online computational strategy for the RB ``compliant'' output. In Section 4, we define our \emph{a posteriori} error estimators for our RB approach, and provide the computation procedures for the two ingredients of our error estimators, which are the dual norm of the residual and the coercivity lower bound. In Section 5, we briefly discuss the extension of our RB methodology to the ``non-compliant'' output. In Section 6, we show several numerical results to illustrate the capability of this method, with a final subsection devoted to provide an introduction to more complex nonlinear problems. Finally, in Section 7, we draw discussions and news on future works.

\section{Preliminaries}
\label{sec:2}
In this Section we shall first present a ``unified'' formulation for all the linear elasticity cases -- for isotropic and orthotropic materials, 2D Cartesian and 3D axisymmetric configurations -- we consider in this study. We then introduce a domain decomposition and geometric mapping strategy to recast the formulation in the ``affine forms'', which is a crucial requirement for our RB approximation. Finally, we define the ``truth'' finite element approximation, upon which we shall build the RB approximation, introduced in the next Section.

\subsection{Formulation on the ``Original'' Domain}
\label{subsec:2}

\subsubsection{Isotropic/Orthotropic materials}
We first briefly describe our problem formulation based on the original settings (denoted by a superscript $^{\rm o}$). We consider a solid body in two dimensions $\Omega^{\rm o}(\bfmu) \in \mathbb{R}^2$ with boundary $\Gamma^{\rm o}$, where $\bfmu \in \calD \subset \mathbb{R}^P$ is the input parameter {and} $\calD$ is the parameter domain \cite{sneddon1999mathematical,sneddon2000mathematical}. For the sake of simplicity, in this section, we assume implicitly that any ``original'' quantities (stress, strain, domains, boundaries, etc) with superscript $^{\rm o}$ will depend on the input parameter $\bfmu$, e.g. $\Omega^{\rm o} \equiv \Omega^{\rm o}(\bfmu)$. 

We first make the following assumptions: (i) the solid is free of body forces, (ii) there {are} negligible thermal strains; note that the extension to include either or both body forces/thermal strains is straightforward. Let us denote $u^{\rm o}$ as the displacement field, and the spatial coordinate $\bfx^{\rm o} = (x^{\rm o}_1,x^{\rm o}_2)$, the linear elasticity equilibrium reads
\begin{equation}\label{eqn:equlibrium}
\frac{\partial \sigma^{\rm o}_{ij}}{\partial x^{\rm o}_j} = 0, \quad \text{in} \ \Omega^{\rm o}
\end{equation}
where $\sigma^{\rm o}$ denotes the stresses, which are related to the strains $\varepsilon^{\rm o}$
by
\begin{equation*}
    \sigma^{\rm o}_{ij} = C_{ijkl}\varepsilon^{\rm o}_{kl},  \quad 1 \leq i,j,k,l \leq 2
\end{equation*}
where
\begin{equation*}
\varepsilon^{\rm o}_{kl} = \dfrac{1}{2}\bigg(\frac{\partial u^{\rm o}_k}{\partial x^{\rm o}_l} + \frac{\partial u^{\rm o}_l}{\partial x^{\rm o}_k} \bigg),
\end{equation*}
$u^{\rm o} = (u^{\rm o}_1,u^{\rm o}_2)$ is the displacement and $C_{ijkl}$ is the elastic tensor, which can be expressed in a matrix form as
\begin{equation*}
    [\bfC] = \left[
              \begin{array}{cccc}
                C_{1111} & C_{1112} & C_{1121} & C_{1122} \\
                C_{1211} & C_{1212} & C_{1221} & C_{1222} \\
                C_{2111} & C_{2112} & C_{2121} & C_{2122} \\
                C_{2211} & C_{2212} & C_{2221} & C_{2222} \\
              \end{array}
            \right]
    = [\bfB]^T [\bfE] [\bfB],
\end{equation*}
where
\begin{equation*}
    [\bfB] = \left[
              \begin{array}{cccc}
                1 & 0 & 0 & 0 \\
                0 & 0 & 0 & 1 \\
                0 & 1 & 1 & 0 \\
              \end{array}
            \right] \quad
    [\bfE] = \left[
              \begin{array}{ccc}
                c_{11} & c_{12} & 0 \\
                c_{21} & c_{21} & 0 \\
                0 & 0 & c_{33} \\
              \end{array}
            \right].
\end{equation*}
The matrix $[\bfE]$ varies for different material types and is given in the Appendix.

We next consider Dirichlet boundary conditions for both components of $u^{\rm o}$:
\begin{equation*}
u^{\rm o}_i = 0 \quad {\rm{on}} \quad \Gamma_{D,i}^{\rm o},
\end{equation*}
and Neumann boundary conditions:
\begin{eqnarray*}
    \sigma^{\rm o}_{ij}e^{\rm o}_{n,j} &=& \bigg\{
                                             \begin{array}{ccc}
                                               f^{\rm o}_ne^{\rm o}_{n,i} & \rm{on} &  \Gamma^{\rm o}_{N} \\
                                               0 & \rm{on} & \Gamma^{\rm o}\backslash\Gamma^{\rm o}_{N}\\
                                             \end{array}
\end{eqnarray*}
where $f^{\rm o}_n$ is the specified stress on boundary edge $\Gamma^{\rm o}_{N}$ respectively; and $\bfe^{\rm o}_n = [e^{\rm o}_{n,1}, e^{\rm o}_{n,2}]$ is the unit normal on $\Gamma^{\rm o}_{N}$. Zero value of $f^{\rm o}_n$ indicate free stress (homogeneous Neumann conditions) on a specific boundary. Here we only consider homogeneous Dirichlet boundary conditions, but extensions to non-homogeneous Dirichlet boundary conditions and/or nonzero traction Neumann boundary conditions are simple and straightforward.

We then introduce the functional space
\begin{equation*}
    X^{\rm o} = \{v = (v_1,v_2) \in (H^1(\Omega^{\rm o}))^2 \ | \ v_i = 0 \ {\rm on} \ \Gamma^{\rm o}_{D,i} , i = 1,2\},
\end{equation*}
here $H^1(\Omega^{\rm o}) = \{ v \in L^2(\Omega^{\rm o}) \ | \ \nabla v \in (L^2(\Omega^{\rm o}))^2 \}$ {and} $L^2(\Omega^{\rm o})$ is the space of square-integrable functions over $\Omega^{\rm o}$. By multiplying \refeq{eqn:equlibrium} by a test function $v \in X^{\rm o}$ and integrating by part over $\Omega^{\rm o}$ we obtain the weak form
\begin{equation}\label{eqn:wf0}
    \int_{\Omega^{\rm o}}\frac{\partial v_i}{\partial x^{\rm o}_j}{C}_{ijkl}\frac{\partial
    u^{\rm o}_k}{\partial x^{\rm o}_l}d\Omega^{\rm o} =
    \int_{\Gamma^{\rm o}_{N}}{f}^{\rm o}_n e^{\rm o}_{n,j} v_jd\Gamma^{\rm o}.
\end{equation}

Finally, we define our output of interest, which usually is a measurement (of our displacement field or even equivalent derived solutions such as stresses, strains) over a boundary segment $\Gamma^{\rm o}_L$ or a part of the domain $\Omega^{\rm o}_L$. Here we just consider a simple case,
\begin{equation}\label{eqn:output0}
    s^{\rm o}(\bfmu) = \int_{\Gamma^{\rm o}_L}f^{\rm o}_{\ell,i}u^{\rm o}_id\Gamma^{\rm o},
\end{equation}
i.e the measure of the displacement on either or both $x^{\rm o}_1$ and $x^{\rm o}_2$ direction along $\Gamma^{\rm o}_L$ with multipliers $f^{\rm o}_{\ell,i}$; more general forms for the output of interest can be extended straightforward. Note that our output of interest is a linear function of the displacement; extension to quadratic function outputs can be found in \cite{huynh07:ijnme}.

We can then now recover our abstract statement: Given a $\bfmu \in \calD$, we evaluate
$$s^{\rm o}(\bfmu) = \ell^{\rm o}(u^{\rm o};\bfmu),$$
where $u^{\rm o} \in X^{\rm o}$ satisfies
$$a^{\rm o}(u^{\rm o},v;\bfmu) = f^{\rm o}(v;\bfmu), \quad \forall v \in X^{\rm o}.$$
Here $a^{\rm o}(w,v;\bfmu): X^{\rm o} \times X^{\rm o} \rightarrow \mathbb{R}$, $\forall w,v \in X^{\rm o}$ is the symmetric and positive bilinear form associated to the left hand side term of \refeq{eqn:wf0}; $f^{\rm o}(v;\bfmu): X^{\rm o} \rightarrow \mathbb{R}$ and $\ell^{\rm o}(v;\bfmu): X^{\rm o} \rightarrow \mathbb{R}$, $\forall v \in X^{\rm o}$ are the linear forms associated to the right hand side terms of \refeq{eqn:wf0} and \refeq{eqn:output0}, respectively. It shall be proven convenience to recast $a^{\rm o}(\cdot,\cdot;\bfmu)$, $f^{\rm o}(\cdot;\bfmu)$ and $\ell^{\rm o}(\cdot;\bfmu)$ in the following forms
\begin{equation}\label{eqn:bilinear_o}
    a^{\rm o}(w,v;\bfmu) = \int_{\Omega^{\rm o}}
    \left[\dfrac{\partial w_1}{\partial x^{\rm o}_1}, \dfrac{\partial w_1}{\partial x^{\rm o}_2}, \dfrac{\partial w_2}{\partial x^{\rm o}_1}, \dfrac{\partial w_2}{\partial x^{\rm o}_2}, w_1 \right]
    [\bfS^a]
    \left[
      \begin{array}{c}
        \dfrac{\partial v_1}{\partial x^{\rm o}_1} \\ \dfrac{\partial v_1}{\partial x^{\rm o}_2} \\
        \dfrac{\partial v_2}{\partial x^{\rm o}_1} \\ \dfrac{\partial v_2}{\partial x^{\rm o}_2} \\
        v_1 \\
      \end{array}
    \right]
    d\Omega^{\rm o}, \ \forall w,v \in X^{\rm o},
\end{equation}
\begin{equation}\label{eqn:linear_fo}
    f^{\rm o}(v;\bfmu) =
    \int_{\Gamma^{\rm o}_N}
    [\bfS^f]
    \left[
      \begin{array}{c}
        v_1 \\ v_2 \\
      \end{array}
    \right]
    d\Gamma^{\rm o}, \quad \forall v \in X^{\rm o},
\end{equation}
\begin{equation}\label{eqn:linear_lo}
    \ell^{\rm o}(v;\bfmu) = \int_{\Gamma^{\rm o}_L}
    [\bfS^{\ell}]
    \left[
      \begin{array}{c}
        v_1 \\ v_2 \\
      \end{array}
    \right]
    d\Gamma^{\rm o}, \quad \forall v \in X^{\rm o},
\end{equation}
where $[\bfS^a] \in \mathbb{R}^{5\times 5}$; $[\bfS^f]\in \mathbb{R}^2$ and $[\bfS^{\ell}]\in \mathbb{R}^2$ are defined as
\begin{equation*}
    [\bfS^a] = \left[
               \begin{array}{ll}
                 [\bfC] & [\boldsymbol{0}]^{4\times 1} \\
                 \left[\boldsymbol{0}\right]^{1\times 4} & 0 \\
               \end{array}
             \right], \quad
    [\bfS^{f}] = \left[
               \begin{array}{cc}
                 f^{\rm o}_ne^{\rm o}_{n,1} & f^{\rm o}_ne^{\rm o}_{n,2} \\
               \end{array}
             \right], \quad
    [\bfS^{\ell}] = \left[
               \begin{array}{cc}
                 f^{\rm o}_{\ell,1}& f^{\rm o}_{\ell,2} \\
               \end{array}
             \right].
\end{equation*}

\subsubsection{Axisymmetric}
Now we shall present the problem formulation for the axisymmetric case. In a cylindrical coordinate system $(r,z,\theta),$\footnote{For the sake of simple illustration, we omit the ``original'' superscript $^{\rm o}$ on $(r,z,\theta)$.} the elasticity equilibrium reads
\begin{eqnarray*}
    \dfrac{\partial \sigma^{\rm o}_{rr}}{\partial r} + \dfrac{\partial \sigma^{\rm o}_{zr}}{\partial z} + \dfrac{\sigma^{\rm o}_{rr}-\sigma^{\rm o}_{\theta\theta}}{r} &=& 0 , \quad \text{in} \quad \Omega^{\rm o} \\
    \dfrac{\partial \sigma^{\rm o}_{rz}}{\partial r} + \dfrac{\partial \sigma^{\rm o}_{zz}}{\partial z} + \dfrac{\sigma^{\rm o}_{rz}}{r} &=& 0 \nonumber, \quad \text{in} \quad \Omega^{\rm o}
\end{eqnarray*}
where $\sigma^{\rm o}_{rr}$, $\sigma^{\rm o}_{zz}$, $\sigma^{\rm o}_{rz}$, $\sigma^{\rm o}_{\theta\theta}$ are the stress components given by
\begin{equation*}
    \left[
      \begin{array}{c}
        \sigma^{\rm o}_{rr} \\
        \sigma^{\rm o}_{zz} \\
        \sigma^{\rm o}_{\theta\theta} \\
        \sigma^{\rm o}_{rz} \\
      \end{array}
    \right] = \underbrace{\dfrac{E}{(1+\nu)(1-2\nu)}\left[\begin{array}{cccc}
                (1-\nu) & \nu & \nu & 0 \\
                \nu & (1-\nu) & \nu & 0 \\
                \nu & \nu & (1-\nu) & 0 \\
                0 & 0 & 0 & \dfrac{1-2\nu}{2} \\
              \end{array}\right]}_{[\bfE]}
    \left[
      \begin{array}{c}
        \varepsilon^{\rm o}_{rr} \\
        \varepsilon^{\rm o}_{zz} \\
        \varepsilon^{\rm o}_{\theta\theta} \\
        \varepsilon^{\rm o}_{rz} \\
      \end{array}
    \right]\nonumber,
\end{equation*}
where $E$ and $\nu$ are the axial Young's modulus and Poisson ratio, respectively. We only consider isotropic material, however, extension to general to anisotropic material is possible; as well as axisymmetric plane stress and plane strain \cite{zienkiewics05:FEM1}.
The strain $\varepsilon^{\rm o}_{rr}$, $\varepsilon^{\rm o}_{zz}$, $\varepsilon^{\rm o}_{rz}$, $\varepsilon^{\rm o}_{\theta\theta}$ are given by
\begin{equation}\label{eqn:stress_axs}
    \left[
      \begin{array}{c}
        \varepsilon^{\rm o}_{rr} \\
        \varepsilon^{\rm o}_{zz} \\
        \varepsilon^{\rm o}_{\theta\theta} \\
        \varepsilon^{\rm o}_{rz} \\
      \end{array}
    \right] = \left[
      \begin{array}{c}
        \dfrac{\partial u^{\rm o}_r}{\partial r} \\[6pt]
        \dfrac{\partial u^{\rm o}_z}{\partial z} \\[6pt]
        \dfrac{u^{\rm o}_r}{r} \\[6pt]
        \dfrac{\partial u^{\rm o}_r}{\partial z} + \dfrac{\partial u^{\rm o}_z}{\partial r} \\
      \end{array}
    \right],
\end{equation}
where $u^{\rm o}_r$, $u^{\rm o}_z$ are the radial displacement and axial displacement, respectively.

Assuming that the axial axis is $x^{\rm o}_2$, let $[u^{\rm o}_1,u^{\rm o}_2] \equiv [\dfrac{u^{\rm o}_r}{r}, u^{\rm o}_z]$ and denoting \\ $[x^{\rm o}_1, x^{\rm o}_2, x^{\rm o}_3] \equiv [r,z,\theta]$, we can then express \refeq{eqn:stress_axs} as
\begin{equation*}
    \left[
      \begin{array}{c}
        \varepsilon^{\rm o}_{11} \\
        \varepsilon^{\rm o}_{22} \\
        \varepsilon^{\rm o}_{33} \\
        \varepsilon^{\rm o}_{12} \\
      \end{array}
    \right] = [\hat{\bfE}]\underbrace{\left[\begin{array}{ccccc}
                x^{\rm o}_1 & 0 & 0 & 0 & 1 \\
                0 & 0 & 0 & 1 & 0 \\
                0 & 0 & 0 & 0 & 1 \\
                0 & x^{\rm o}_1 & 1 & 0 & 0 \\
              \end{array}\right]}_{[\bfB_a]}
     \left[
      \begin{array}{c}
        \dfrac{\partial u^{\rm o}_1}{\partial x^{\rm o}_1} \\[7pt]
        \dfrac{\partial u^{\rm o}_1}{\partial x^{\rm o}_2} \\[7pt]
        \dfrac{\partial u^{\rm o}_2}{\partial x^{\rm o}_1} \\[7pt]
        \dfrac{\partial u^{\rm o}_2}{\partial x^{\rm o}_2} \\[7pt]
        u^{\rm o}_1 \\
      \end{array}
    \right]\nonumber.
\end{equation*}

As in the previous case, we consider the usual homogeneous Dirichlet boundary conditions on $\Gamma^{\rm o}_{D,i}$ and Neumann boundary conditions on $\Gamma^{\rm o}$.
Then if we consider the output of interest $s^{\rm o}(\bfmu)$ defined upon $\Gamma^{\rm o}_L$, we arrive at the same abstract statement where
\begin{equation*}
[\bfS^{a}] = x^{\rm o}_1[\bfB_a]^T[\bfE][\bfB_a], \
[\bfS^{f}] = \left[(x^{\rm o}_1)^2f^{\rm o}_ne^{\rm o}_{n,1} , \ x^{\rm o}_1f^{\rm o}_ne^{\rm o}_{n,2} \right], \
[\bfS^{\ell}] = \left[ x^{\rm o}_1f^{\rm o}_ne^{\rm o}_{n,1} ,\  f^{\rm o}_ne^{\rm o}_{n,2}\right].
\end{equation*}

Note that the $x^{\rm o}_1$ multipliers appear in $[\bfS^{f}]$ during the weak form derivation, while in $[\bfS^{\ell}]$, in order to retrieve the measurement for the axial displacement $u^{\rm o}_r$ rather than $u^{\rm o}_1$ due to the change of variables. Also, the $2\pi$ multipliers in both $a^{\rm o}(\cdot,\cdot;\bfmu)$ and $f^{\rm o}(\cdot;\bfmu)$ are disappeared in the weak form during the derivation, and can be included in $\ell^{\rm o}(\cdot;\bfmu)$, i.e. incorporated to $[\bfS^{\ell}]$ if measurement is required to be done in thruth (rather than in the axisymmetric) domain.

\subsection{Formulation on Reference Domain}

The RB requires that the computational domain must be parameter-independent; however, our ``original'' domain $\Omega^{\rm o}(\bfmu)$ is obviously parameter-dependent. Hence, to transform $\Omega^{\rm o}(\bfmu)$ into the computational domain, or ``reference'' (parameter-independent) domain $\Omega$, we must perform geometric transformations in order to express the bilinear and linear forms in our abstract statement in appropriate ``affine forms''. This ``affine forms'' formulation allows us to model all possible configurations, corresponding to every $\bfmu \in \calD$, based on a single reference-domain \cite{Quarteroni2011,rozza08:ARCME}.

\subsubsection{Geometry Mappings}

We first assume that, for all $\bfmu \in \calD$, $\Omega^{\rm o}(\bfmu)$ is expressed as
\begin{equation*}
    \Omega^{\rm o}(\bfmu) = \bigcup_{s = 1}^{L_{\rm reg}}\Omega^{\rm o}_s(\bfmu),
\end{equation*}
where the $\Omega^{\rm o}_s(\bfmu)$, $s = 1,\ldots,L_{\rm reg}$ are mutually non-overlapping subdomains. In two dimensions, $\Omega^{\rm o}_s(\bfmu)$, $s = 1,\ldots,L_{\rm reg}$ is a set of triangles (or in the general case, a set of ``curvy triangles''\footnote{In fact, a ``curvy triangle'' \cite{rozza08:ARCME} is served as the building block. For its implementation see} \cite{rbMIT_URL}.) such that all important domains/edges (those defining different material regions, boundaries, pressures/tractions loaded boundary segments, or boundaries which the output of interests are calculated upon) are included in the set. In practice, such a set is generated by a constrained Delaunay triangulation.

We next assume that there exists a reference domain $\Omega \equiv \Omega^{\rm o}(\bfmu_{\rm ref}) =
\bigcup_{s=1}^{L_{\rm reg}}\Omega_s$ where, for any $\bfx^{\rm o} \in \Omega_s$, $s = 1,\ldots,L_{\rm reg}$, its image $\bfx^{\rm o} \in \Omega^{\rm o}_s(\bfmu)$ is given by
\begin{equation}\label{eqn:local_mapping}
    \bfx^{\rm o}(\bfmu)
    = \mathcal{T}^{\rm aff}_s(\bfmu;\bfx)
    = [\bfR^{\rm aff}_s(\bfmu)]\bfx + [\bfG^{\rm aff}_s(\bfmu)],
\end{equation}
where $[\bfR^{\rm aff}_s(\bfmu)] \in \mathbb{R}^{2 \times 2}$ and $[\bfG^{\rm aff}_s(\bfmu)] \in \mathbb{R}^2$. It thus follows from our definitions that $\mathcal{T}_s(\bfmu;\bfx):\Omega_s\rightarrow\Omega_s^{\rm o}$, $1 \leq s \leq L_{\rm reg}$ is an (invertible) affine mapping from $\Omega_s$ to $\Omega^{\rm o}_s(\bfmu)$, hence the Jacobian $|{\rm det}([\bfR^{\rm aff}_s(\bfmu)])|$ is strictly positive, and that the derivative transformation matrix, $[\bfD^{\rm aff}_s(\bfmu)] = [\bfR^{\rm aff}_s(\bfmu)]^{-1}$ is well defined. We thus can write
\begin{equation}\label{eqn:spatial_derv}
    \frac{\partial}{\partial x^{\rm o}_i} = \frac{\partial x_j}{\partial x^{\rm o}_i}\frac{\partial}{\partial x_j}
    = D^{\rm aff}_{s,ij}(\bfmu)\frac{\partial}{\partial x_j}, \quad 1 \leq i,j \leq 2.
\end{equation}
As in two dimensions, an affine transformation maps a triangle to a triangle, we can readily calculate $[\bfR^{\rm aff}_s(\bfmu)]$ and $[\bfG^{\rm aff}_s(\bfmu)]$ for each subdomains $s$ by simply solving a systems of six equations forming from \refeq{eqn:local_mapping} by matching parametrized coordinates to reference coordinates for the three triangle vertices.

We further require a mapping continuity condition: for all $\bfmu \in \mathcal{D}$,
\begin{equation*}
    \mathcal{T}_s(\bfmu;\bfx) = \mathcal{T}_{s'}(\bfmu;\bfx),
    \quad \forall \bfx \in \Omega_s \cap \Omega_{s'}, \quad 1 \leq
    s,s' \leq L_{\rm reg}.
\end{equation*}
This condition is automatically held if there is no curved edge in the set of $\Omega^{\rm o}_s(\bfmu)$. If a domain contains one or more ``important'' curved edge, special ``curvy triangles'' must be generated appropriately to honour the continuity condition. We refer the readers to \cite{rozza08:ARCME} for the full discussion and detail algorithm for such cases.

The global transformation is for $\bfx \in \Omega$, the image $\bfx^{\rm o} \in \Omega^{\rm o}(\bfmu)$ is given by
\begin{equation*}
    \bfx^{\rm o}(\bfmu) = \mathcal{T}(\bfmu;\bfx).
\end{equation*}
It thus follows that $\mathcal{T}(\bfmu;\bfx):\Omega\rightarrow\Omega^{\rm o}(\bfmu)$ is a piecewise-affine geometric mapping.

\subsubsection{Affine Forms}

We now define our functional space $X$ as 
\begin{equation*}
    X = \{v = (v_1,v_2) \in (H^1(\Omega))^2 | v_i = 0 \ {\rm on} \ \Gamma_{D,i}, i = 1,2\},
\end{equation*}
and recast our bilinear form $a^{\rm o}(w,v;\bfmu)$, by invoking \refeq{eqn:bilinear_o}, \refeq{eqn:local_mapping} and \refeq{eqn:spatial_derv} to obtain $\forall w,v \in X(\Omega)$
\begin{eqnarray*}
a(w,v;\bfmu) &=& \int_{\bigcup_{s=1}^{L_{\rm reg}}\Omega_s}\left[
        \dfrac{\partial w_1}{\partial x_1},
        \dfrac{\partial w_1}{\partial x_2},
        \dfrac{\partial w_2}{\partial x_1},
        \dfrac{\partial w_2}{\partial x_2},
        w_1
    \right][\bfS^{a,\rm aff}_s(\bfmu)]
    \left[
      \begin{array}{c}
        \dfrac{\partial v_1}{\partial x_1} \\
        \dfrac{\partial v_1}{\partial x_2} \\
        \dfrac{\partial v_2}{\partial x_1} \\
        \dfrac{\partial v_2}{\partial x_2} \\
        v_1 \\
      \end{array}
    \right]d\Omega. 
\end{eqnarray*}
where $[\bfS^{a,\rm aff}_s(\bfmu)] = [\bfH_s(\bfmu)][\bfS^a_s][\bfH_s(\bfmu)]^T|{\rm det}([\bfR^{\rm aff}_s(\bfmu)])|$ is the effective elastic tensor matrix, in which
\begin{equation*}
    [\bfH_s(\bfmu)] = \left(
                        \begin{array}{ccc}
                          [\bfD_s(\bfmu)] & [\boldsymbol{0}]^{2 \times 2} & 0 \\
                          \phantom{1}[\boldsymbol{0}]^{2 \times 2} & [\bfD_s(\bfmu)] & 0 \\
                          0 & 0 & 1 \\
                        \end{array}
                      \right).
\end{equation*}
Similarly, the linear form $f^{\rm o}(v;\bfmu)$, $\forall v \in X$ can be transformed as
\begin{eqnarray*}
f(v;\bfmu) = \int_{\bigcup_{s=1}^{L_{\rm reg}}\Gamma_{N_s}}[\bfS^{f,{\rm aff}}_s]
    \left[
      \begin{array}{c}
        v_1 \\ v_2 \\
      \end{array}
    \right]
    d\Gamma,
\end{eqnarray*}
where $\Gamma_{N_s}$ denotes the partial boundary segment of $\Gamma_N$ of the subdomain $\Omega_s$ and $[\bfS^{f,{\rm aff}}_s] = \|([\bfR_s(\bfmu)]\bfe_n)\|_2[\bfS^f]$ is the effective load vector, where $\bfe_n$ is the normal vector to $\Gamma_{N_s}$ and $\|\cdot\|_2$ denotes the usual Euclidean norm. The linear form $\ell(v;\bfmu)$ is also transformed in the same manner.

We then replace all ``original'' $x^{\rm o}_1$ and $x^{\rm o}_2$ in the effective elastic tensor matrix $[\bfS_s^{a,\rm aff}(\bfmu)]$, effective load/output vectors $[\bfS^{f,\rm aff}_s(\bfmu)]$ and $[\bfS^{\ell,\rm aff}_s(\bfmu)]$ by \refeq{eqn:local_mapping} to obtain a $\bfx^{\rm o}$-free effective elastic tensor matrix and effective load/output vectors, respectively.\footnote{Here we note that, the Young's modulus $E$ in the isotropic and axisymmetric cases (or $E_1$, $E_2$ and $E_3$ in the orthotropic case only} in certain conditions) can be a polynomial function of the spatial coordinates $\bfx^{\rm o}$ as well, and we still be able to obtain our affine forms \refeq{eqn:affine}.

We next expand the bilinear form $a(w,v;\bfmu)$ by treating each entry of the effective elastic tensor matrix for each subdomain separately, namely
\begin{eqnarray}\label{eqn:bilinear_exp}
    a(w,v;\bfmu) &=& S_{1,11}^{a,\rm aff}(\bfmu)\int_{\Omega_1}\frac{\partial w_1}{\partial x_1}\frac{\partial v_1}{\partial x_1} +  S_{1,12}^{a,\rm aff}(\bfmu)\int_{\Omega_1}\frac{\partial w_1}{\partial x_1}\frac{\partial v_1}{\partial x_2} + \ldots \\
     && + S_{L_{\rm reg},55}^{a,\rm aff}(\bfmu)\int_{\Omega_{L_{\rm reg}}}w_1w_1.
\end{eqnarray}
Note that here for simplicity, we consider the case where there is no spatial coordinates in $[\bfS^{\ell,\rm aff}_s(\bfmu)]$. In general (especially for axisymmetric case), some or most of the integrals may take the form of $\int_{\Omega_s}(x_1)^m(x_2)^n\dfrac{\partial w_i}{\partial x_j}\dfrac{\partial v_k}{\partial x_l}$, where $m,n \in \mathbb{R}$.

Taking into account the symmetry of the bilinear form and the effective elastic tensor matrix, there will be at most $Q^a = 7L_{\rm reg}$ terms in the expansion. However, in practice, most of the terms can be collapsed by noticing that not only there will be a lot of zero entries in $[\bfS_s^{a,\rm aff}(\bfmu)]$, $s = 1,\ldots,L_{\rm reg}$, but also there will be a lot of duplicated or ``linearly dependent'' entries, for example, $S_{1,11}^{a,\rm aff}(\bfmu) = [{\rm Const}]S_{2,11}^{a,\rm aff}(\bfmu)$. We can then apply a symbolic manipulation technique \cite{rozza08:ARCME} to identify, eliminate all zero terms in \refeq{eqn:bilinear_exp} and collapse all ``linear dependent'' terms to end up with a minimal $Q^a$ expansion. The same procedure is also applied for the linear forms $f(\cdot;\bfmu)$ and $\ell(\cdot;\bfmu)$.

Hence the abstract formulation of the linear elasticity problem in the reference domain $\Omega$ reads as follow: given $\bfmu \in \calD$, find
\begin{equation*}
    s(\bfmu) = \ell(u(\bfmu);\bfmu),
\end{equation*}
where $u(\bfmu) \in X$ satisfies
\begin{equation*}
    a(u(\bfmu),v;\bfmu) = f(v;\bfmu), \quad \forall v \in X,
\end{equation*}
where all the bilinear and linear forms are in affine forms,
\begin{eqnarray}\label{eqn:affine}
    a(w,v;\bfmu) &=& \sum_{q=1}^{Q^a}\Theta_q^a(\bfmu) a_q(w,v), \nonumber \\
    f(v;\bfmu) &=& \sum_{q=1}^{Q^f}\Theta_q^f(\bfmu) f_q(v), \nonumber \\
    \ell(v;\bfmu) &=& \sum_{q=1}^{Q^{\ell}}\Theta_q^{\ell}(\bfmu) \ell_q(v), \quad \forall w,v, \in X.
\end{eqnarray}
Here $\Theta_q^a(\bfmu)$, $a_q(w,v)$, $q = 1,\ldots,Q^a$, $f_q(v)$; $\Theta_q^f(\bfmu)$, $f_q(v)$, $q = 1,\ldots,Q^f$, and $\Theta_q^{\ell}(\bfmu)$, $\ell_q(v)$, $q = 1,\ldots,Q^{\ell}$ are parameter-dependent coefficient and parameter-independent bilinear and linear forms, respectively.

We close this section by defining several useful terms. We first define our inner product and energy norm as
\begin{equation}\label{eqn:inner_prod}
(w,v)_X = a(w,v;\overline{\bfmu})
\end{equation}
and $\|w\|_X = (w,w)^{1/2}$, $\forall w,v \in X$, respectively, where $\overline{\bfmu} \in \calD$ is an arbitrary parameter. Certain other inner norms and associated norms are also possible \cite{rozza08:ARCME}. We then define our coercivity and continuity constants as
\begin{equation}\label{eqn:inf}
    \alpha(\bfmu) = \inf_{w\in X}\frac{a(w,v;\bfmu)}{\|w\|_X^2},
\end{equation}
\begin{equation}\label{eqn:sup}
    \gamma(\bfmu) = \sup_{w\in X}\frac{a(w,v;\bfmu)}{\|w\|_X^2},
\end{equation}
respectively. We assume that $a(\cdot,\cdot;\bfmu)$ is symmetric, $a(w,v;\bfmu) = a(v,w;\bfmu)$, $\forall w,v \in X$, coercive, $\alpha(\bfmu) > \alpha_0 > 0$, and continuous, $\gamma(\bfmu) < \gamma_0 < \infty$; and also our $f(\cdot;\bfmu)$ and $\ell(\cdot;\bfmu)$ are bounded functionals. It follows that problem which is well-defined and has a unique solution. Those conditions are automatically satisfied given the nature of our considered problems \cite{sneddon1999mathematical, sneddon2000mathematical}.

\subsection{Truth approximation}

From now on, we shall restrict our attention to the ``compliance'' case ($f(\cdot;\bfmu) = \ell(\cdot;\bfmu)$). Extension to the non-compliance case will be discuss in the Section~5.

We now apply the finite element method and we provide a matrix formulation \cite{CMCS-CONF-2009-002}: given $\bfmu \in \calD$, we evaluate
\begin{equation}\label{eqn:FE_out}
    s(\bfmu) = [\bfF^\calN(\bfmu)]^T[\bfu^\calN(\bfmu)],
\end{equation}
where $[\bfu^\calN(\bfmu)]$ represents a finite element solution $u^{\calN}(\bfmu) \in X^{\calN} \in X$ of size $\calN$ which satisfies
\begin{equation}\label{eqn:FE_stiff}
    [\bfK^\calN(\bfmu)][\bfu^\calN(\bfmu)] = [\bfF^\calN(\bfmu)];
\end{equation}
here $[\bfK^\calN(\bfmu)]$, and $[\bfF^\calN(\bfmu)]$ and the (discrete forms) stiffness matrix and load vector of $a(\cdot,\cdot;\bfmu)$, and $f(\cdot;\bfmu)$, respectively. Note that the stiffness matrix $[\bfK^\calN(\bfmu)]$ is symmetric positive definite (SPD). By invoking the affine forms \refeq{eqn:affine}, we can express $[\bfK^\calN(\bfmu)]$, and $[\bfF^\calN(\bfmu)]$ as
\begin{eqnarray}\label{eqn:affine_FE}
    [\bfK^\calN(\bfmu)] &=& \sum_{q = 1}^{Q^a}\Theta_q^a(\bfmu)[\bfK^\calN_q], \nonumber \\
    \left[\bfF^\calN(\bfmu)\right] &=& \sum_{q = 1}^{Q^f}\Theta_q^f(\bfmu)[\bfF^\calN_q],
\end{eqnarray}
where $[\bfK^\calN_q]$, $[\bfF^\calN_q]$ and are the discrete forms of the parameter-independent bilinear and linear forms $a_q(\cdot,\cdot)$ and $f_q(\cdot)$, respectively. We also denote (the SPD matrix) $[\bfY^\calN]$ as the discrete form of our inner product \refeq{eqn:inner_prod}. We also assume that the size of of our FE approximation, $\calN$ is large enough such that our FE solution is an accurate approximation of the exact solution.

\section{Reduced Basis Method}

In this Section we shall restrict our attention by recalling the RB  method for the ``compliant'' output. We shall first define the RB spaces and the Galerkin projection. We then describe an Offline-Online computational strategy, which allows us to obtain $\calN$-independent calculation of the RB output approximation \cite{hesthaven2015certified,NgocCuong2005}.

\subsection{RB Spaces and Greedy algorithm}

To define the RB approximation we first introduce a (nested) Lagrangian parameter sample for $1 \leq N \leq N_{\max}$,
\begin{equation*}
    S_N = \{\bfmu_1,\bfmu_2,\ldots,\bfmu_N\},
\end{equation*}
and associated hierarchical reduced basis spaces $(X_N^\calN =) W^\calN_N$, $1 \leq N \leq N_{\max}$,
\begin{equation*}
    W^\calN_N = {\rm span}\{u^{\calN}(\bfmu_n),1 \leq n \leq N\} ,
\end{equation*}
where $\bfmu_n \in \calD$ are determined by the means of a Greedy sampling algorithm \cite{rozza08:ARCME, quarteroni2015reduced}; this is an interarive procedure where at each step a new basis function is added in order to improve the precision of the basis set.

 The key point of this methodology is the availability of an estimate of the error induced by replacing the full space $X^{\calN}$ with the reduced order one $W^\calN_N$ in the variational formulation. More specifically we assume that for all $\bfmu \in \calD$ there exist an estimator $\eta(\bfmu)$ such that
\begin{equation*}
|| u^{\calN}(\bfmu) - u^\calN_{{\rm RB},N}(\bfmu)|| \leq \eta(\bfmu) , 
\end{equation*} 
where $u^{\calN}(\bfmu) \in X^{\calN} \in X$ represents the finite element solution, $u^\calN_{{\rm RB},N}(\bfmu)\in X_N^\calN \subset X^\calN$ the reduced basis one and we can choose either the induced or the energy norm.

During this iterative basis selection process and if at the j-th step a j-dimensional reduced basis space $W^\calN_j$ is given, the next basis function is the one that maximizes the estimated model order reduction error given the j-dimensional space $W^\calN_j$  over $\calD$. So at the $n+1$ iteration we select $$\bfmu_{n+1} = arg \max_{\bfmu \in \calD} \eta(\bfmu)$$ and compute $u^{\calN}(\bfmu_{n+1})$ to enrich the reduced space. This is repeated until the maximal estimated error is below a required error tolerance. With this choice the greedy algorithm always selects the next parameter sample point as the one for which the model error is the maximum as estimated by $\eta(\bfmu)$ and this yields a basis that aims to be optimal in the maximum norm over $\calD$.

Furthermore we can rewrite the reduced space as 
\begin{equation*}
W^\calN_N = {\rm span}\{\zeta^{\calN}_n,1 \leq n \leq N\},
\end{equation*}
where the basis functions $\left\{\zeta^{\calN}\right\}$ are computed from the snapshots $u^{\calN}(\bfmu)$ by a Gram-Schmidt orthonormalization process such that $[\bfzeta^{\calN}_m]^T[\bfY^\calN][\bfzeta^{\calN}_n] = \delta_{mn}$, where $\delta_{mn}$ is the Kronecker-delta symbol. We then define our orthonormalized-snapshot matrix $[\bfZ_N] \equiv [\bfZ^\calN_N] = [[\bfzeta^{\calN}_1]|\cdots|[\bfzeta^{\calN}_n]]$ of dimension $\calN \times N$.

\subsection{Galerkin Projection}

We then apply a Galerkin projection on our ``truth'' problem \cite{almroth78:_autom,noor81:_recen,noor82,noor80:_reduc,rozza08:ARCME}: given $\bfmu \in \calD$, we could evaluate the RB output as
\begin{equation*}
    s_N(\bfmu) = [\bfF^\calN(\bfmu)]^T[\bfu^\calN_{{\rm RB},N}(\bfmu)],
\end{equation*}
where
\begin{equation}\label{eqn:RB_sol}
[\bfu^\calN_{{\rm RB},N}(\bfmu)] = [\bfZ_N][\bfu_N(\bfmu)]
\end{equation}
represents the RB solution $\bfu^\calN_{{\rm RB},N}(\bfmu) \in X_N^\calN \subset X^\calN$ of size $\calN$. Here $[\bfu_N(\bfmu)]$ is the RB coefficient vector of dimension $N$ satisfies the RB ``stiffness'' equations
\begin{equation}\label{eqn:RB_semifull}
    [\bfK_N(\bfmu)][\bfu_N(\bfmu)] = [\bfF_N(\bfmu)],
\end{equation}
where
\begin{eqnarray}\label{eqn:RB_comp1}
[\bfK_N(\bfmu)] &=& [\bfZ_N]^T[\bfK^\calN(\bfmu)][\bfZ_N], \nonumber \\
\left[\bfF_N(\bfmu)\right] &=& [\bfZ_N]^T[\bfF^\calN(\bfmu)].
\end{eqnarray}
Note that the system \refeq{eqn:RB_semifull} is of small size: it is just a set of $N$ linear algebraic equations, in this way we can now evaluate our output as
\begin{equation}\label{eqn:RB_outsemifull}
    s_N(\bfmu) = [\bfF_N(\bfmu)]^T[\bfu_N(\bfmu)].
\end{equation}
It can be shown \cite{patera07:book} that the condition number of the RB ``stiffness'' matrix $[\bfZ_N]^T[\bfK^\calN(\bfmu)][\bfZ_N]$ is bounded by $\gamma_0(\bfmu)/\alpha_0(\bfmu)$, and independent of both $N$ and $\calN$.

\subsection{Offline-Online Procedure}

Although the system \refeq{eqn:RB_semifull} is of small size, the computational cost for assembling the RB ``stiffness'' matrix (and the RB ``output'' vector $[\bfF^\calN(\bfmu)]^T[\bfZ_N]$) is still involves $\calN$ and costly, $O(N\calN^2 + N^2\calN)$ (and $O(N\calN)$, respectively). However, we can use our affine forms \refeq{eqn:affine} to construct very efficient Offline-Online procedures, as we shall discuss below.

We first insert our affine forms \refeq{eqn:affine_FE} into the expansion \refeq{eqn:RB_semifull} and \refeq{eqn:RB_outsemifull}, by using \refeq{eqn:RB_comp1} we obtain
\begin{equation*}
    \sum_{q = 1}^{Q^a}\Theta_q^a(\bfmu)[\bfK_{qN}][\bfu_N(\bfmu)]
    = \sum_{q = 1}^{Q^f}\Theta_q^f(\bfmu)[\bfF_{qN}]
\end{equation*}
and
\begin{equation*}
    s_N(\bfmu) = \sum_{q = 1}^{Q^f}\Theta_q^f(\bfmu)[\bfF_{qN}][\bfu_N(\bfmu)],
\end{equation*}
respectively. Here
\begin{eqnarray*}
    [\bfK_{qN}] &=& [\bfZ_N]^T[\bfK^\calN_q][\bfZ_N], \quad 1 \leq q \leq Q^a \\
    \left[\bfF_{qN}\right] &=& [\bfZ_N]^T[\bfF^\calN_q], \quad 1 \leq q \leq Q^f,
\end{eqnarray*}
are parameter independent quantities that can be computed just once and than stored for all the subsequent $\bfmu$-dependent queries.
We then observe that all the ``expensive'' matrices $[\bfK_{qN}]$, $1 \leq q \leq Q^a$, $1 \leq N \leq N_{\max}$ and vectors $[\bfF_{qN}]$, $1 \leq q \leq Q^f$, $1 \leq N \leq N_{\max}$, are now separated and parameter-independent, hence those can be \emph{pre-computed} in an Offline-Online procedure.

In the Offline stage, we first compute the $[\bfu^\mathcal{N}(\bfmu^n)]$, $1 \leq n \leq N_{\max}$, form the matrix $[\bfZ_{N_{\max}}]$ and then form and store $[\bfF_{N_{\max}}]$ and $[\bfK_{qN_{\max}}]$. The Offline operation count depends on $N_{\max}$, $Q^a$ and $\mathcal{N}$ but requires only $O(Q^aN_{\max}^2 + Q^fN_{\max} + Q^\ell N_{\max})$ permanent storage.

In the Online stage, for a given $\bfmu$ and $N$ ($1 \leq N \leq N_{\max}$), we retrieve the pre-computed $[\bfK_{qN}]$ and $[\bfF_{N}]$ (subarrays of $[\bfK_{qN_{\max}}]$, $[\bfF_{N_{\max}}]$), form $[\bfK_N(\bfmu)]$, solve the resulting $N \times N$ system \refeq{eqn:RB_semifull} to obtain $\{\bfu_N(\bfmu)\}$, and finally evaluate the output $s_N(\bfmu)$ from \refeq{eqn:RB_outsemifull}. The Online operation count is thus $O(N^3)$ and independent of $\mathcal{N}$. The implication of the latter is two-fold: first, we will achieve very fast response in the many-query and real-time contexts, as $N$ is typically very small, $N \ll \mathcal{N}$; and second, we can choose $\mathcal{N}$ arbitrary large -- to obtain as accurate FE predictions as we wish -- without adversely affecting the Online (marginal) cost.

\section{\emph{A posteriori} error estimation}

In this Section we recall the \emph{a posteriori} error estimator for our RB approximation. We shall discuss in details the computation procedures for the two ingredients of the error estimator: the dual norm of the residual and the coercivity lower bound. We first present the Offline-Online strategy for the computation of the dual norm of the residual; we then briefly discuss the Successive Constraint Method \cite{huynh07:cras} in order to compute the coercivity lower bound.

\subsection{Definitions}

We first introduce the error $e^\calN(\bfmu) \equiv u^\calN(\bfmu) - u^\calN_{{\rm RB},N}(\bfmu) \in X^\calN$ and the
residual $r^\calN(v;\bfmu) \in (X^\calN)'$ (the dual space to $X^\calN)$, $\forall v \in X^\calN$,
\begin{equation}\label{eqn:residual}
    r^\calN(v;\bfmu) = f(v) - a(u^\calN(\bfmu),v;\bfmu),
\end{equation}
which can be given in the discrete form as
\begin{equation}\label{eqn:FE_residual}
    [\bfr^\calN(\bfmu)] = [\bfF^\calN(\bfmu)] - [\bfK^\calN(\bfmu)][\bfu^\calN_{{\rm RB},N}(\bfmu)].
\end{equation}
We then introduce the Riesz representation of $r^\calN(v;\bfmu)$: $\hat{e}(\bfmu) \in X^\calN$ defined by $(\hat{e}(\bfmu),v)_{X^\calN} = r^\calN(v;\bfmu)$, $\forall v \in X^\calN$. In vector form, $\hat{e}(\bfmu)$ can be expressed as
\begin{equation}\label{eqn:rres}
    [\bfY^\calN][\bfehat(\bfmu)] = [\bfr^{\calN}(\bfmu)].
\end{equation}

We also require a lower bound to the coercivity constant
\begin{equation}\label{eqn:inf_FE}
\alpha^\calN(\bfmu) = \inf_{w \in X^{\calN}}\frac{a(w,w;\bfmu)}{\|w\|^2_{X^{\calN}}},
\end{equation}
such that $0< \alpha_{\rm LB}^\calN(\bfmu) \leq \alpha^\calN(\bfmu)$, $\forall \bfmu \in \calD$.

We may now define our error estimator for our output as
\begin{equation}
    \Delta_N^s(\bfmu) \equiv \frac{\|\hat{e}(\bfmu)\|^2_{X^\calN}}{\alpha^\calN_{\rm LB}},
\end{equation}
where $\|\hat{e}(\bfmu)\|_{X^\calN}$ is the dual norm of the residual.
We can also equip the error estimator with an effectivity defined by
\begin{equation}
    \eta_N^s(\bfmu) \equiv \frac{\Delta_N^s(\bfmu)}{|s^\calN(\bfmu)-s_N(\bfmu)|}.
\end{equation}
We can readily demonstrate \cite{rozza08:ARCME, patera07:book} that
\begin{equation*}
    1 \leq \eta_N^s(\bfmu) \leq \frac{\gamma_0(\bfmu)}{\alpha^\calN_{\rm LB}(\bfmu)}, \quad \forall \bfmu \in \calD;
\end{equation*}
so that the error estimator is both \emph{rigorous} and \emph{sharp}. Note that here we can only claim the \emph{sharp} property for this current ``compliant'' case.

We shall next provide procedures for the computation of the two ingredients of our error estimator: we shall first discuss the Offline-Online strategy to compute the dual norm of the residual $\|\hat{e}(\bfmu)\|_{X^\calN}$, and then provide the construction for the lower bound of the coercivity constant $\alpha^\calN(\bfmu)$.

\subsection{Dual norm of the residual}

In discrete form, the dual norm of the residual $\varepsilon(\bfmu) = \|\hat{e}(\bfmu)\|_{X^\calN}$ is given by
\begin{equation}\label{eqn:dnres}
    \varepsilon^2(\bfmu) = [\bfehat(\bfmu)]^T[\bfY^\calN][\bfehat(\bfmu)].
\end{equation}
We next invoke \refeq{eqn:FE_residual}, \refeq{eqn:rres} and \refeq{eqn:dnres} to arrive at
\begin{eqnarray}\label{eqn:dnres_der1}
    \varepsilon^2(\bfmu) &=& \bigg([\bfF^\calN(\bfmu)] - [\bfK^\calN(\bfmu)][\bfu^\calN_{{\rm RB},N}(\bfmu)])\bigg)^T[\bfY^\calN]^{-1}\bigg([\bfF^\calN(\bfmu)] - [\bfK^\calN(\bfmu)][\bfu^\calN_{{\rm RB},N}(\bfmu)]\bigg) \nonumber \\
    &=& [\bfF^\calN(\bfmu)]^T[\bfY^\calN]^{-1}[\bfF^\calN(\bfmu)] - 2[\bfF^\calN(\bfmu)]^T[\bfY^\calN]^{-1}[\bfK^\calN(\bfmu)] \nonumber  \\
    && + [\bfK^\calN(\bfmu)]^T[\bfY^\calN]^{-1}[\bfK^\calN(\bfmu)].
\end{eqnarray}
We next defines the ``pseudo''-solutions $[\bfP^f_{q}] = [\bfY^\calN]^{-1}[\bfF_{q}^\calN]$, $1 \leq q \leq Q^f$ and $[\bfP^a_{qN}] = [\bfY^\calN]^{-1}[\bfK_q^\calN][\bfZ_N]$, $1 \leq q \leq Q^a$, then apply the affine form \refeq{eqn:affine_FE} and \refeq{eqn:RB_sol} into \refeq{eqn:dnres_der1} to obtain
\begin{eqnarray}\label{eqn:dnres_der2}
    \varepsilon^2(\bfmu) &=& \sum_{q=1}^{Q^f}\sum_{q'=1}^{Q^f}\Theta_q^f(\bfmu)\Theta_{q'}^f(\bfmu)\bigg([\bfP^f_{q}]^T[\bfY^\calN][\bfP^f_{q'}]\bigg) \\ && -2\sum_{q=1}^{Q^a}\sum_{q'=1}^{Q^f}\Theta_q^a(\bfmu)\Theta_{q'}^f(\bfmu)\bigg([\bfP^f_{q}]^T[\bfY^\calN][\bfP^a_{q'N}]\bigg)[\bfu_N^{RB}(\bfmu)] \nonumber \\
    &&+\sum_{q=1}^{Q^a}\sum_{q'=1}^{Q^a}\Theta_q^a(\bfmu)\Theta_{q'}^a(\bfmu)[\bfu_N^{RB}(\bfmu)]^T\bigg([\bfP^a_{qN}]^T[\bfY^\calN][\bfP^a_{q'N}]\bigg)[\bfu_N^{RB}(\bfmu)] \nonumber.
\end{eqnarray}
It is observed that all the terms in bracket in \refeq{eqn:dnres_der2} are all parameter-independent, hence they can be \emph{pre-computed} in the Offline stage. The Offline-Online strategy is now clear.

In the Offline stage we form the parameter-independent quantities. We first compute the ``pseudo''-solutions $[\bfP^f_{q}] = [\bfY^\calN]^{-1}[\bfF_{q}^\calN]$, $1 \leq q \leq Q^f$ and $[\bfP^a_{qN}] = [\bfY^\calN]^{-1}[\bfK_q^\calN][\bfZ_N]$, $1 \leq q \leq Q^a$, $1 \leq N \leq N_{\max}$; and form/store $[\bfP^f_{q}]^T[\bfY^\calN][\bfP^f_{q'}]$, $1 \leq q, q' \leq Q^f$, $[\bfP^f_{q}]^T[\bfY^\calN][\bfP^a_{q'N}]$, $1 \leq q \leq Q^f$, $1 \leq q \leq Q^a$, $1 \leq N \leq N_{\max}$,\\ $[\bfP^a_{qN}][\bfY^\calN][\bfP^a_{q'N}]$, $1 \leq q, q' \leq Q^a$, $1 \leq N \leq N_{\max}$. The Offline operation count depends on $N_{\max}$, $Q^a$, $Q^f$, and $\calN$.

In the Online stage, for a given $\bfmu$ and $N$ ($1 \leq N \leq N_{\max}$), we retrieve the pre-computed quantities $[\bfP^f_{q}]^T[\bfY^\calN][\bfP^f_{q'}]$, $1 \leq q, q' \leq Q^f$, $[\bfP^f_{q}]^T[\bfY^\calN][\bfP^a_{q'N}]$, $1 \leq q \leq Q^f$, $1 \leq q \leq Q^a$, and  $[\bfP^a_{qN}]^T[\bfY^\calN][\bfP^a_{q'N}]$, $1 \leq q, q' \leq Q^a$, and then evaluate the sum \refeq{eqn:dnres_der2}. The Online operation count is dominated by $O(((Q^a)^2+(Q^f)^2)N^2)$ and independent of $\mathcal{N}$.

\subsection{Lower bound of the coercivity constant}

We now briefly address some elements for the computation of the lower bound in the coercive case. In order to derive the discrete form of the coercivity constant $\refeq{eqn:inf_FE}$ we introduce the discrete eigenvalue problem: given $\bfmu \in \calD$, find the minimum set $([\bfchi_{\min}(\bfmu)],\lambda_{\min}(\bfmu))$ such that
\begin{eqnarray}\label{eqn:inf_truth}
    [\bfK^\calN(\bfmu)][\bfchi(\bfmu)] &=& \lambda_{\min}[\bfY^\calN][\bfchi(\bfmu)], \nonumber \\
    \left[\bfchi(\bfmu)\right]^T[\bfY^\calN][\bfchi(\bfmu)] &=& 1.
\end{eqnarray}
We can then recover
\begin{equation}
    \alpha^\calN(\bfmu) = \sqrt{\lambda_{\min}(\bfmu)}.
\end{equation}
However, the eigenproblem $\refeq{eqn:inf_truth}$ is of size $\calN$, so using direct solution as an ingredient for our error estimator is very expensive. Hence, we will construct an inexpensive yet of good quality lower bound $\alpha_{\rm LB}^\calN(\bfmu)$ and use this lower bound instead of the truth (direct) expensive coercivity constant $\alpha^\calN(\bfmu)$ in our error estimator.

For our current target problems, our bilinear form is coercive and symmetric. We shall construct our coercivity lower bound by the Successive Constraint Method (SCM) \cite{huynh07:cras}. It is noted that the SCM method can be readily extended to non-symmetric as well as non-coercive bilinear forms \cite{huynh07:cras,rozza08:ARCME,patera07:book,huynh08:infsupLB}.

We first introduce an alternative (albeit not very computation-friendly) discrete form for our coercivity constant as
\begin{eqnarray}\label{eqn:inf_Y}
    {\rm minimum}  && \sum_{q = 1}^{Q^a}\Theta_q^a(\bfmu)y_q, \\
    {\rm subject \ to} && y_q = \frac{[\bfw_q]^T[\bfK^\calN_q][\bfw_q]}{[\bfw_q]^T[\bfY^\calN][\bfw_q]}, \quad 1 \leq q \leq Q^a, \nonumber
\end{eqnarray}
where $[\bfw_q]$ is the discrete vector of any arbitrary $w_y \in X^\calN$.

We shall now ``relax'' the constraint in \refeq{eqn:inf_Y} by defining the ``continuity constraint box'' associated with $y_{q,\min}$ and $y_{q,\max}$, $1 \leq q \leq Q^a$ obtained from the minimum set $([\bfy_-(\bfmu)],y_{q,\min})$ and maximum set $([\bfy_+(\bfmu)],y_{q,\max})$ solutions of the eigenproblems
\begin{eqnarray*}
    [\bfK^\calN_q][\bfy_-(\bfmu)] &=& y_{q,\min}[\bfY^\calN][\bfy_-(\bfmu)], \\
    \left[\bfy_-(\bfmu)\right][\bfY^\calN][\bfy_-(\bfmu)] &=& 1,
\end{eqnarray*}
and
\begin{eqnarray*}
    [\bfK^\calN_q][\bfy_+(\bfmu)] &=& y_{q,\max}[\bfY^\calN][\bfy_+(\bfmu)], \\
    \left[\bfy_+(\bfmu)\right][\bfY^\calN][\bfy_+(\bfmu)] &=& 1,
\end{eqnarray*}
respectively, for $1 \leq q \leq Q^a$. We next define a ``coercivity constraint'' sample
\begin{equation*}
    C_J = \{\bfmu^{\rm SCM}_1 \in \calD, \ldots, \bfmu^{\rm SCM}_J \in \calD\},
\end{equation*}
and denote $C_J^{M,\bfmu}$ the set of $M$ $(1 \leq M \leq J)$ points in $C_J$ closest (in the usual Euclidean norm) to a given $\bfmu \in \calD$. The construction of the set $C_J$ is done by means of a Greedy procedure \cite{huynh07:cras,rozza08:ARCME,patera07:book}. The Greedy selection of $C_J$ can be called the ``Offline stage'', which involves the solutions of $J$ eigenproblems \refeq{eqn:inf_truth} to obtain $\alpha^\calN(\bfmu)$, $\forall \bfmu \in C_J$.

We may now define our lower bound $\alpha^\calN_{\rm LB}(\bfmu)$ as the solution of
\begin{eqnarray}\label{eqn:inf_LB}
    {\rm minimum}  && \sum_{q = 1}^{Q^a}\Theta_q^a(\bfmu)y_q, \\
    {\rm subject \ to} && y_{q,\min} \leq y_q \leq y_{q,\max}, \quad 1 \leq q \leq Q^a, \nonumber \\
    && \sum_{q = 1}^{Q^a}\Theta_q^a(\bfmu')y_q \geq \alpha^\calN(\bfmu'), \quad \forall \bfmu' \in C_J^{M,\bfmu}. \nonumber
\end{eqnarray}

We then ``restrict'' the constraint in \refeq{eqn:inf_Y} and define our upper bound $\alpha^\calN_{\rm UB}(\bfmu)$ as the solution of
\begin{eqnarray}\label{eqn:inf_UB}
    {\rm mininum}  && \sum_{q = 1}^{Q^a}\Theta_q^a(\bfmu)y_{q,*}(\bfmu'), \\
    {\rm subject \ to} && y_{q,*}(\bfmu') = [\bfchi(\bfmu')]^T[\bfK^\calN_q][\bfchi(\bfmu')], \quad 1 \leq q \leq Q^a, \quad \forall \bfmu' \in C_J^{M,\bfmu}, \nonumber
\end{eqnarray}
where $[\bfchi(\bfmu)]$ is defined by \refeq{eqn:inf_truth}. It can be shown \cite{huynh07:cras,rozza08:ARCME,patera07:book} that the feasible region of \refeq{eqn:inf_UB} is a subset of that of \refeq{eqn:inf_Y}, which in turn, is a subset of that of \refeq{eqn:inf_LB}: hence $\alpha^\calN_{\rm LB}(\bfmu) \leq \alpha^\calN(\bfmu) \leq \alpha^\calN_{\rm UB}(\bfmu)$.

We note that the lower bound \refeq{eqn:inf_LB} is a linear optimization problem (or Linear Program (LP)) which contains $Q^a$ design variables and $2Q^a + M$ inequality constraints. Given a value of the parameter $\bfmu$, the Online evaluation $\bfmu \rightarrow \alpha^\calN_{\rm LB}(\bfmu)$ is thus as follows: we find the subset $C_J^{M,\bfmu}$ of $C_J$ for a given $M$, we then calculate $\alpha^\calN_{\rm LB}(\bfmu)$ by solving the LP \refeq{eqn:inf_LB}. The crucial point here is that the online evaluation $\bfmu \rightarrow \alpha^\calN_{\rm LB}(\bfmu)$ is totally independent of $\calN$. The upper bound \refeq{eqn:inf_LB}, however, can be obtained as the solution of just a simple enumeration problem; the online evaluation of $\alpha^\calN_{\rm UB}(\bfmu)$ is also independent of $\calN$. In general, the upper bound $\alpha^\calN_{\rm UB}(\bfmu)$ is not used in the calculation of the error estimator, however, it is used in the Greedy construction of the set $C_J$ \cite{huynh07:cras}. In practice, when the set $C_J$ does not guarantee to produce a positive $\alpha^\calN_{\rm LB}(\bfmu)$, the upper bound $\alpha^\calN_{\rm UB}(\bfmu)$ can be used as a substitution for $\alpha^\calN_{\rm UB}(\bfmu)$ since it approximates the ``truth'' $\alpha^\calN(\bfmu)$ in a very way; however we will lose the rigorous property of the error estimators.

\section{Extension of the RB method to non-compliant output}

We shall briefly provide the extension of our RB methodology for the ``non-compliant'' case in this Section. We first present a suitable primal-dual formulation for the ``non-compliant'' output; we then briefly provide the extension to the RB methodology, including the RB approximation and its \emph{a posteriori} error estimation.

\subsection{Adjoint Problem}

We shall briefly discuss the extension of our methodology to the non-compliant problems. We still require that both $f$ and $\ell$ are bounded functionals, but now $(f(\cdot;\bfmu) \neq \ell(\cdot;\bfmu))$. We still use the previous abstract statement in Section~2. We begin with the definition of the dual problem associated to $\ell$: find $\psi(\bfmu) \in X$ (our ``adjoint'' or ``dual'' field) such that
\begin{equation*}
    a(v,\psi(\bfmu);\bfmu) = -\ell(\bfmu), \quad \forall v \in X.
\end{equation*}

\subsection{Truth approximation}

We now again apply the finite element method to the dual formulation: given $\bfmu \in \calD$, we evaluate
\begin{equation*}
    s(\bfmu) = [\bfL^\calN(\bfmu)]^T[\bfu^\calN(\bfmu)],
\end{equation*}
where $[\bfu^\calN(\bfmu)]$ is the finite element solution of size $\calN$ satisfying \refeq{eqn:FE_stiff}. The discrete form of the dual solution $\psi^\calN(\bfmu) \in X^\calN$ is given
\begin{equation*}
    [\bfK^\calN(\bfmu)][\bfpsi^\calN(\bfmu)] = -[\bfL^\calN(\bfmu)];
\end{equation*}
here $[\bfL^\calN(\bfmu)]$ is the discrete load vector of $\ell(\cdot;\bfmu)$. We also invoke the affine forms \refeq{eqn:affine} to express $[\bfL^\calN(\bfmu)]$ as
\begin{eqnarray}\label{eqn:affine_FE_out}
    [\bfL^\calN(\bfmu)] &=& \sum_{q = 1}^{Q^\ell}\Theta_q^\ell(\bfmu)[\bfL^\calN_q],
\end{eqnarray}
where all the $[\bfL^\calN_q]$ are the discrete forms of the parameter-independent linear forms $\ell_q(\cdot)$, $1 \leq q \leq Q^\ell$.

\subsection{Reduced Basis Approximation}

We now define our RB spaces: we shall need to define two Lagrangian parameter samples set, $S_{N^{\rm pr}} = \{\bfmu_1,\bfmu_2,\ldots,\bfmu_{N^{\rm pr}}\}$ and $S_{N^{\rm du}}= \{\bfmu_1,\bfmu_2,\ldots,\bfmu_{N^{\rm du}}\}$ corresponding to the set of our primal and dual parameter samples set, respectively. We also associate the primal and dual reduced basis spaces $(X_{N^{\rm pr}}^\calN =) W^\calN_{N^{\rm pr}}$, $1 \leq N \leq N^{\rm pr}_{\max}$ and $(X_{N^{\rm du}}^\calN =) W^\calN_{N^{\rm du}}$, $1 \leq N \leq N^{\rm du}_{\max}$ to our $S_{N^{\rm pr}}$ and $S_{N^{\rm du}}$ set, respectively, which are constructed from the primal $u^{\calN}(\bfmu)$ and dual $\psi^{\calN}(\bfmu)$ snapshots by a Gram-Schmidt process as in Section~3. Finally, we denote our primal and dual orthonormalized-snapshot as $[\bfZ^{\rm pr}_{N^{\rm pr}}]$ and $[\bfZ^{\rm du}_{N^{\rm du}}]$ basis matrices, respectively.

\subsection{Galerkin Projection}

We first denote the RB primal approximation to the primal ``truth'' approximation $u^\calN(\bfmu)$ as $u_{{\rm RB},N}^\calN(\bfmu)$ and the RB dual approximation to the primal ``truth'' dual approximation $\psi^\calN(\bfmu)$ as $\psi_{{\rm RB},N}^\calN(\bfmu)$: their discrete forms are given by $[\bfu_{RB,N^{\rm pr}}^\calN(\bfmu)] = [\bfZ^{\rm pr}_{N^{\rm pr}}][\bfu_{N^{\rm pr}}(\bfmu)]$ and $[\bfpsi_{RB,N^{\rm du}}^\calN(\bfmu)] = [\bfZ^{\rm du}_{N^{\rm du}}][\bfpsi_{N^{\rm du}}(\bfmu)]$, respectively.

We then apply a Galerkin projection (note that in this case, a Galerkin-Petrov projection is also possible \cite{rozza08:ARCME, patera07:book,benner2015}). given a $\bfmu \in \calD$, we evaluate the RB output
\begin{equation*}
    s_{N^{\rm pr},N^{\rm du}}(\bfmu) = [\bfL^\calN(\bfmu)]^T[\bfu^\calN_{{\rm RB},N^{\rm pr}}(\bfmu)] - [\bfr^\calN_{\rm pr}(\bfmu)]^T[\bfpsi^\calN_{{\rm RB},N^{\rm du}}(\bfmu)],
\end{equation*}
recall that $[\bfr^\calN_{\rm pr}(\bfmu)]$ is the discrete form of the RB primal residual defined in \refeq{eqn:residual}. The RB coefficient primal and dual are given by
\begin{eqnarray}\label{eqn:semifull_du}
    \sum_{q=1}^{Q^a}\Theta_q^a(\bfmu)[\bfK_{qN^{\rm pr}N^{\rm pr}}][\bfu_{N^{\rm pr}}(\bfmu)] &=& \sum_{q=1}^{Q^f}\Theta_q^f(\bfmu)[\bfF_{qN^{\rm pr}}], \nonumber \\
    \sum_{q=1}^{Q^a}\Theta_q^a(\bfmu)[\bfK_{qN^{\rm du}N^{\rm du}}[\bfpsi_{N^{\rm du}}(\bfmu)] &=& -\sum_{q=1}^{Q^\ell}\Theta_q^\ell(\bfmu)[\bfL_{qN^{\rm du}}].
\end{eqnarray}
Note that the two systems \refeq{eqn:semifull_du} are also of small size: their sizes are of $N^{\rm pr}$ and $N^{\rm du}$, respectively. We can now evaluate our output as
\begin{eqnarray}\label{eqn:RB_outsemifull_du}
    s_{N^{\rm pr},N^{\rm du}}(\bfmu) &=& \sum_{q=1}^{Q^\ell}\Theta_q^\ell(\bfmu)[\bfL_{qN^{\rm pr}}][\bfu_{N^{\rm pr}}(\bfmu)] - \sum_{q=1}^{Q^f}\Theta_q^f(\bfmu)[\bfF_{qN^{\rm du}}][\bfpsi_{N^{\rm du}}(\bfmu)] \nonumber\\
    &&+\sum_{q=1}^{Q^a}\Theta_q^a(\bfmu)[\bfpsi_{N^{\rm du}}(\bfmu)]^T[\bfK_{qN^{\rm du}N^{\rm pr}}][\bfu_{N^{\rm pr}}(\bfmu)].
\end{eqnarray}
All the quantities in \refeq{eqn:semifull_du} and \refeq{eqn:RB_outsemifull_du} are given by
\begin{eqnarray*}
[\bfK_{qN^{\rm pr}N^{\rm pr}}] &=& [\bfZ^{\rm pr}_{N^{\rm pr}}]^T[\bfK_q][\bfZ^{\rm pr}_{N^{\rm pr}}], \quad 1 \leq q \leq Q^a, \ 1 \leq N^{\rm pr} \leq N^{\rm pr}_{\max},\\
\left[\bfK_{qN^{\rm du}N^{\rm du}}\right] &=& [\bfZ^{\rm du}_{N^{\rm du}}]^T[\bfK_q][\bfZ^{\rm du}_{N^{\rm du}}], \quad 1 \leq q \leq Q^a, \ 1 \leq N^{\rm du} \leq N^{\rm du}_{\max}, \\
\left[\bfK_{qN^{\rm du}N^{\rm pr}}\right] &=& [\bfZ^{\rm du}_{N^{\rm du}}]^T[\bfK_q][\bfZ^{\rm pr}_{N^{\rm pr}}], \quad 1 \leq q \leq Q^a, \ 1 \leq N^{\rm pr} \leq N^{\rm pr}_{\max}, \ 1 \leq N^{\rm du} \leq N^{\rm du}_{\max} \\
\left[\bfF_{qN^{\rm pr}}\right] &=& [\bfZ^{\rm pr}_{N^{\rm pr}}]^T[\bfF_q], \quad 1 \leq q \leq Q^f, \ 1 \leq N^{\rm pr} \leq N^{\rm pr}_{\max}, \\
\left[\bfF_{qN^{\rm du}}\right] &=& [\bfZ^{\rm du}_{N^{\rm du}}]^T[\bfF_q], \quad 1 \leq q \leq Q^f, \ 1 \leq N^{\rm du} \leq N^{\rm du}_{\max}, \\
\left[\bfL_{qN^{\rm pr}}\right] &=& [\bfZ^{\rm pr}_{N^{\rm pr}}]^T[\bfL_q], \quad 1 \leq q \leq Q^\ell, 1 \leq N^{\rm pr} \leq N^{\rm pr}_{\max}, \\
\left[\bfL_{qN^{\rm du}}\right] &=& [\bfZ^{\rm du}_{N^{\rm du}}]^T[\bfL_q], \quad 1 \leq q \leq Q^\ell, 1 \leq N^{\rm du} \leq N^{\rm du}_{\max}.
\end{eqnarray*}
The computation of the output $s_{N^{\rm pr},N^{\rm du}}(\bfmu)$ clearly admits an Offline-Online computational strategy similar to the one we discuss previously in Section~3.

\subsection{\emph{A posteriori} error estimation}

We now introduce the dual residual $r_{\rm du}^\calN(v;\bfmu)$,
\begin{equation*}
    r_{\rm du}^\calN(v;\bfmu) = -\ell(v) - a(v,\psi_{N^{\rm du}}^\calN(\bfmu);\bfmu), \quad \forall v \in X^\calN.
\end{equation*}
and its Riesz representation of $r_{\rm du}^\calN(v;\bfmu)$: $\hat{e}^{\rm du}(\bfmu) \in X^\calN$ defined by $(\hat{e}^{\rm du}(\bfmu),v)_{X^\calN} = r^\calN_{\rm du}(v;\bfmu)$, $\forall v \in X^\calN$.

We may now define our error estimator for our output as
\begin{equation}
    \Delta_{N^{\rm pr}N^{\rm du}}^s(\bfmu) \equiv \frac{\|\hat{e}^{\rm pr}(\bfmu)\|_{X^\calN}}{(\alpha^\calN_{\rm LB})^{1/2}}\frac{\|\hat{e}^{\rm du}(\bfmu)\|_{X^\calN}}{(\alpha^\calN_{\rm LB})^{1/2}},
\end{equation}
where $\hat{e}^{\rm pr}(\bfmu)$ is the Riesz representation of the primal residual. We then define the effectivity associated with our error bound
\begin{equation}
    \eta_{N^{\rm pr}N^{\rm du}}^s(\bfmu) \equiv \frac{\Delta_{N^{\rm pr}N^{\rm du}}^s(\bfmu)}{|s^\calN(\bfmu)-s_{N^{\rm pr}N^{\rm du}}(\bfmu)|}.
\end{equation}
We can readily demonstrate \cite{rozza08:ARCME, patera07:book, grepl04:_reduc_basis_approx_time_depen} that
\begin{equation*}
    1 \leq \eta_{N^{\rm pr}N^{\rm du}}^s(\bfmu), \quad \forall \bfmu \in \calD;
\end{equation*}
note that the error estimator is still \emph{rigorous}, however it is less \emph{sharp} than that in the ``compliant'' case since in this case we could not provide an upper bound to $\eta_{N^{\rm pr}N^{\rm du}}^s(\bfmu)$.

The computation of the dual norm of the primal/dual residual also follows an Offline-Online computation strategy: the dual norm of the primal residual is in fact, the same as in Section~4.2; the same procedure can be applied to compute the dual norm of the dual residual.

\section{Numerical results}

In this sections we shall consider several ``model problems'' to demonstrate the feasibility of our methodology. We note that in all cases, these model problems are presented in non-dimensional form unless stated otherwise. In all problems below, displacement is, in fact, in non-dimensional form $u = {\tilde{u}\tilde{E}}/{\tilde{\sigma}_0}$, where $\tilde{u}$, $\tilde{E}$, $\tilde{\sigma}_0$ are the dimensional displacement, Young's modulus and load strength, respectively, while $E$ and $\sigma_0$ are our non-dimensional Young's modulus and load strength and usually are around unity.

We shall not provide any details for $\Theta_q^a(\bfmu)$, $\Theta_q^f(\bfmu)$ and $\Theta_q^\ell(\bfmu)$ and their associated bilinear and linear forms $a_q(\cdot,\cdot)$, $f_q(\cdot)$ and $\ell_q(\cdot)$ for any of the below examples as they are usually quite complex, due to the complicated structure of the effective elastic tensor and our symbolic manipulation technique. We refer the users to \cite{huynh07:ijnme, patera07:book, veroy03:_phd_thesis, milani08:RB_LE}, in which all the above terms are provided in details for some simple model problems.

In the below, the timing $t_{\rm FE}$ for an evaluation of the FE solution $\bfmu \rightarrow s^\calN(\bfmu)$ is the computation time taken by solving \refeq{eqn:FE_stiff} and evaluating \refeq{eqn:FE_out} by using \refeq{eqn:affine_FE} and \refeq{eqn:affine_FE_out}, in which all the stiffness matrix components, $[\bfK_q]$, $1\leq q\leq Q^a$, load and output vector components, $[\bfF_q]$, $1\leq q\leq Q^f$ and $[\bfL_q]$, $1\leq q\leq Q^\ell$, respectively, are pre-computed and pre-stored. We do not include the computation time of forming those components (or alternatively, calculate the stiffness matrix, load and output vector directly) in $t_{\rm FE}$.

Finally, for the sake of simplicity, we shall denote the number of basis $N$ defined as $N = N^{\rm pr} = N^{\rm du}$ in all of our model problems in this Section.

\subsection{The arc-cantilever beam}

\begin{figure} [htbp]
\centering
\includegraphics[scale=0.5]{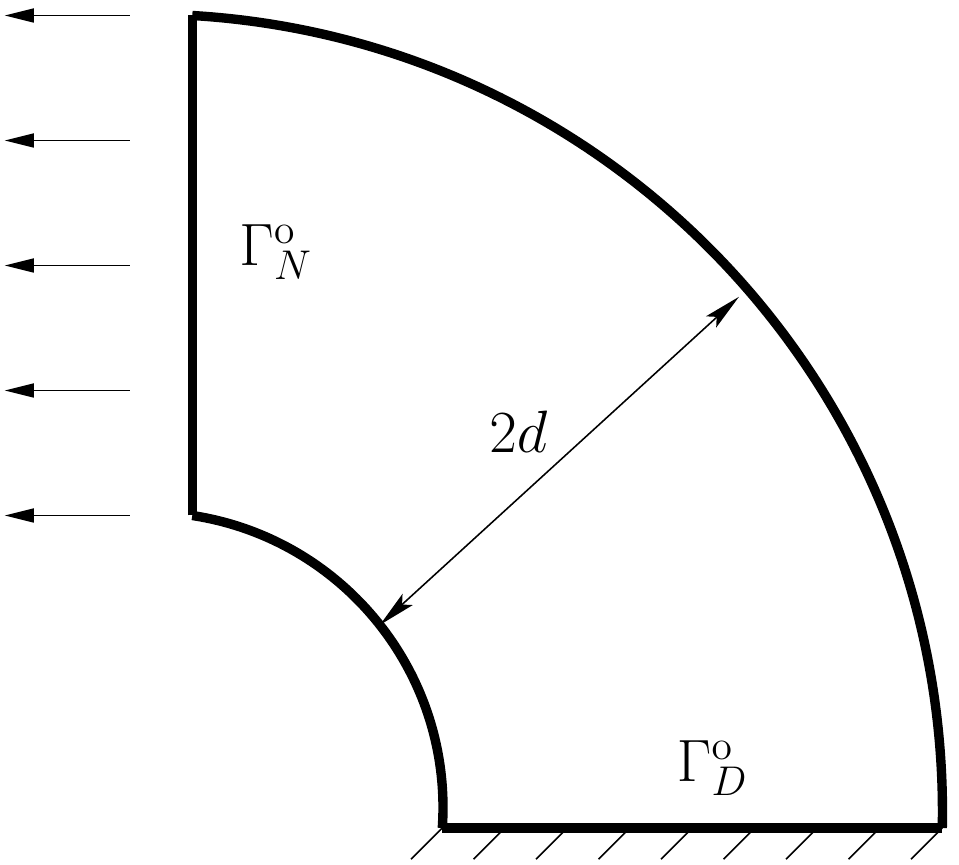}
\caption{The arc-cantilever beam}
\label{fig:ex1_model}
\end{figure}

We consider a thick arc cantilever beam correspond to the domain $\Omega^{\rm o}(\bfmu)$ representing the shape of a quarter of an annulus as shown in Figure~\ref{fig:ex1_model}. We apply (clamped) homogeneous Dirichlet conditions on $\Gamma^{\rm o}_D$ and non-homogeneous Neumann boundary conditions corresponding to a unit tension on $\Gamma^{\rm o}_N$. The width of the cantilever beam is $2d$, and the material is isotropic with $(E,\nu) = (1,0.3)$ under plane stress assumption. Our output of interest is the integral of the tangential displacement ($u_2$) over $\Gamma^{\rm o}_N$, which can be interpreted as the average tangential displacement on $\Gamma^{\rm o}_N$\footnote{The average tangential displacement on $\Gamma^{\rm o}_N$ is not exactly $s(\bfmu)$ but rather $s(\bfmu)/l_{\Gamma^{\rm o}_N}$, where $l_{\Gamma^{\rm o}_N}$ is the length of ${\Gamma^{\rm o}_N}$. It is obviously that the two descriptions of the two outputs, ''integral of'' and ``average of'', are pretty much equivalent to each other.}. Note that our output of interest is ``non-compliant''.

\begin{figure} [htbp]
\centering
\includegraphics[scale=0.3]{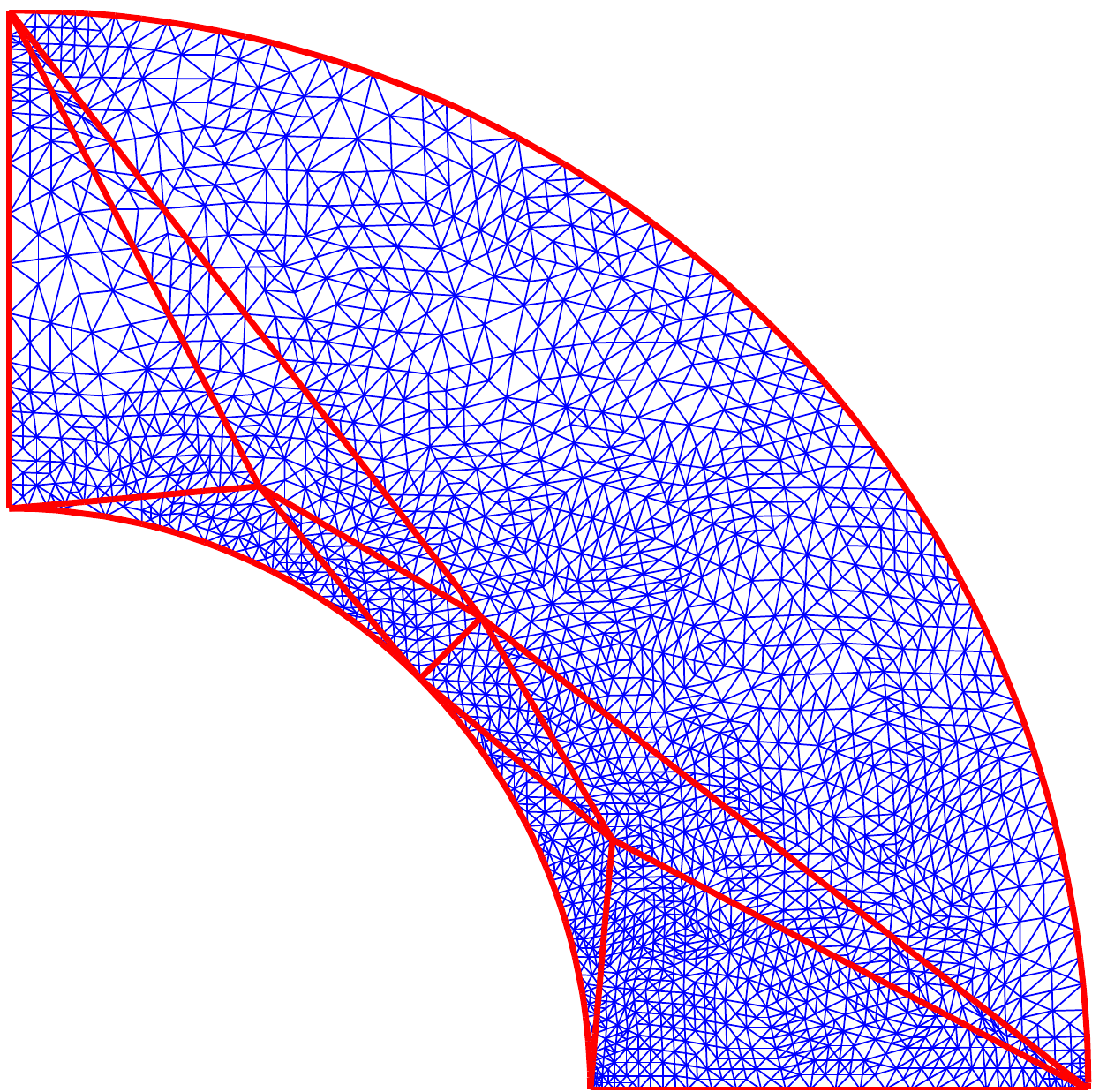}
\caption{The arc-cantilever beam problem: Domain composition and FE mesh}
\label{fig:ex1_mesh}
\end{figure}

The parameter is the half-width of the cantilever beam $\bfmu = [\mu_1] \equiv [d]$. The parameter domain is chosen as $\calD = [0.3, 0.9]$, which can model a moderately thick beam to a very thick beam. We then choose $\bfmu_{\rm ref} = 0.3$ and apply the domain decomposition and obtain $L_{\rm reg} = 9$ subdomains as shown in Figure~\ref{fig:ex1_mesh}, in which three subdomains are the general ``curvy triangles'', generated by our computer automatic procedure \cite{rozza08:ARCME}. Note that geometric transformations are relatively complicated, due to the appearances of the ``curvy triangles'' and all subdomains transformations are classified as the ``general transformation case'' \cite{patera07:book, huynh07:_phd_thesis}. We then recover our affine forms with $Q^a = 54$, $Q^f = 1$ and $Q^l = 1$.

We next consider a FE approximation where the mesh contains $n_{\rm node} = 2747$ nodes and $n_{\rm elem} = 5322$ $P_1$ elements, which corresponds to $\calN = 5426$ degrees of freedoms\footnote{Note that $\calN \neq 2n_{\rm node}$ since Dirichlet boundary nodes are eliminated from the FE system.} as shown in Figure~\ref{fig:ex1_mesh}. To verify our FE approximation, we compare our FE results with the approximated solution for thick arc cantilever beam by Roark \cite{roark01:roark_formula} for a $100$ uniformly distributed test points in $\calD$: the maximum difference between our results and Roark's one is just $2.9\%$.

We then apply our RB approximation. We present in Table~\ref{tab:ex1_tab} our convergence results: the RB error bounds and effectivities as a function of $N (=N^{\rm pr} = N^{\rm du})$. The error bound reported, $\calE_N = \Delta^s_N(\bfmu)/|s_N(\bfmu)|$ is the maximum of the relative error bound over a random test sample $\Xi_{\rm test}$ of size $n_{\rm test} = 100$. We denote by $\overline{\eta}_N^s$ the average of the effectivity $\eta_N^s(\bfmu)$ over $\Xi_{\rm test}$. We observe that average effectivity is of order $O(20-90)$, not very \emph{sharp}, but this is expected due to the fact that the output is ``non-compliant''.
\begin{table}
\centering
\begin{tabular}{|c||c|c|}
    \hline
    $N$ & $\calE_N$ & $\overline{\eta}_N^s$ \\
    \hline
     2 & 3.57\texttt{E}+00 & 86.37 \\
     4 & 3.70\texttt{E}-03 & 18.82 \\
     6 & 4.07\texttt{E}-05 & 35.72 \\
     8 & 6.55\texttt{E}-07 & 41.58 \\
    10 & 1.99\texttt{E}-08 & 40.99 \\
    \hline
\end{tabular}
\caption{The arc-cantilever beam: RB convergence}
\label{tab:ex1_tab}
\end{table}

As regards computational times, a RB online evaluation $\bfmu \rightarrow (s_N(\bfmu),\Delta_N^s(\bfmu))$ requires just $t_{\rm RB} = 115$(ms) for $N = 10$; while the FE solution $\bfmu \rightarrow s^\calN(\bfmu)$ requires $t_{\rm FE} = 9$(s): thus our RB online evaluation is just $1.28\%$ of the FEM computational cost.

\subsection{The center crack problem}

\begin{figure} [htbp]
\centering
\includegraphics[scale=0.5]{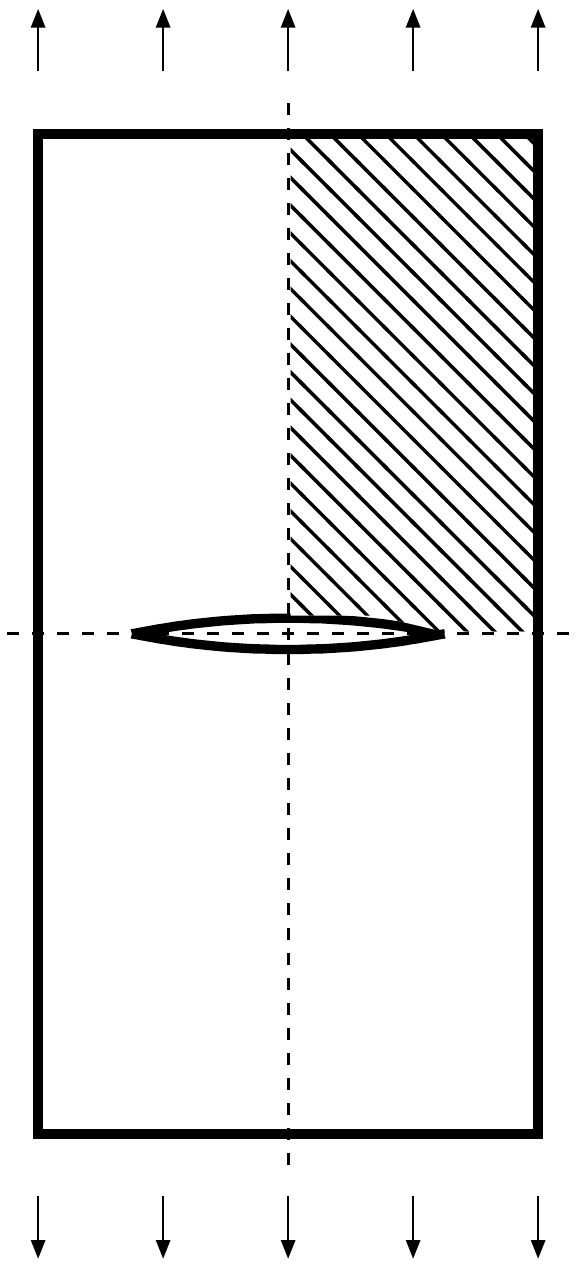}
\caption{The center crack problem}
\label{fig:ex2_modelfull}
\end{figure}

We next consider a fracture model corresponds to a center crack in a plate under tension at both sides as shown in Figure~\ref{fig:ex2_modelfull}.

Due to the symmetry of the geometry and loading, we only consider one quarter of the physical domain, as shown in Figure~\ref{fig:ex2_modelfull}, note that the crack corresponds to the boundary segment $\Gamma^{\rm o}_C$. The crack (in our ``quarter'' model) is of size $d$, and the plate is of height $h$ (and of fixed width $w = 1$). We consider plane strain isotropic material with $(E,\nu) = (1,0.3)$. We consider (symmetric about the $x^{\rm o}_1$ direction and $x^{\rm o}_2$ direction) Dirichlet boundary conditions on the left and bottom boundaries $\Gamma^{\rm o}_L$ and $\Gamma^{\rm o}_B$, respectively; and non-homogeneous Neumann boundary conditions (tension) on the top boundary $\Gamma^{\rm o}_T$. Our ultimate output of interest is the stress intensity factor (SIF) for the crack, which will be derived from an intermediate (compliant) energy output by application of the virtual crack extension approach \cite{parks77:a_stiff_sif}. The SIF plays an important role in the field of fracture mechanics, for examples, if we have to estimate the propagation path of cracks in structures \cite{hutchingson79:fracture}. We further note that analytical result for SIF of a center-crack in a plate under tension is only available for the infinite plate \cite{murakami01:SIFhandbook}, which can be compared with our solutions for small crack length $d$ and large plate height $h$ values.
\begin{figure} [htbp]
\centering
\includegraphics[scale=0.5]{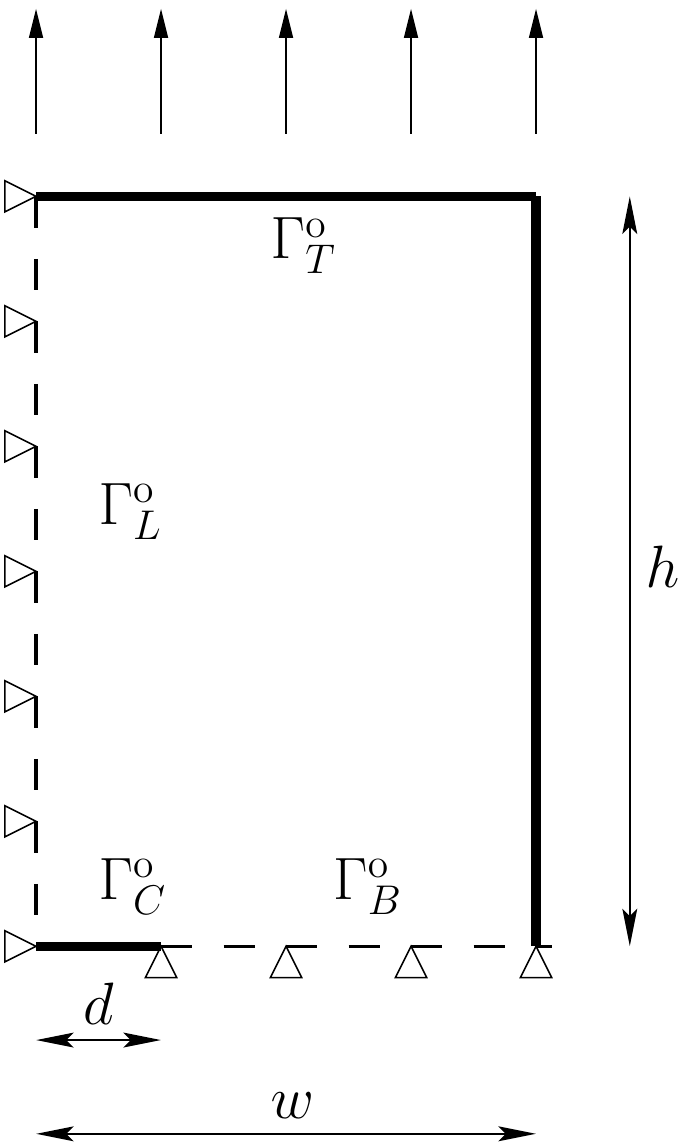}
\caption{The center crack problem}
\label{fig:ex2_model}
\end{figure}

Our parameters are the crack length and the plate height $\bfmu = [\mu_1,\mu_2] \equiv [d, h]$, and the parameter domain is given by $\calD = [0.3,0.7] \times [0.5,2.0]$. We then choose $\bfmu_{\rm ref} = [0.5,1.0]$ and apply a domain decomposition: the final setting contains $L_{\rm reg} = 3$ subdomains, which in turn gives us $Q^a = 10$ and $Q^f = 1$. Note that our ``compliant'' output $s(\bfmu)$ is just an intermediate result for the calculation of the SIF. In particular, the virtual crack extension method (VCE) \cite{parks77:a_stiff_sif} allows us to extract the ``Mode-I'' SIF though the energy $s(\bfmu)$ though the Energy Release Rate (ERR), $G(\bfmu)$, defined by
\begin{equation*}
    G(\bfmu) = -\bigg(\frac{\partial s(\bfmu)}{\partial \mu_1}\bigg).
\end{equation*}
In practice, the ERR is approximated by a finite-difference (FD) approach for a suitable small value $\delta\mu_1$ as
\begin{equation*}
    \widehat{G}(\bfmu) = -\bigg(\frac{s(\bfmu+\delta\mu_1)-s(\bfmu)}{\delta\mu_1}\bigg),
\end{equation*}
which then give the SIF approximation $\widehat{\rm SIF}(\bfmu) = \sqrt{\widehat{G}(\bfmu)/(1-\nu^2)}$.

\begin{figure} [htbp]
\centering
\includegraphics[scale=0.3]{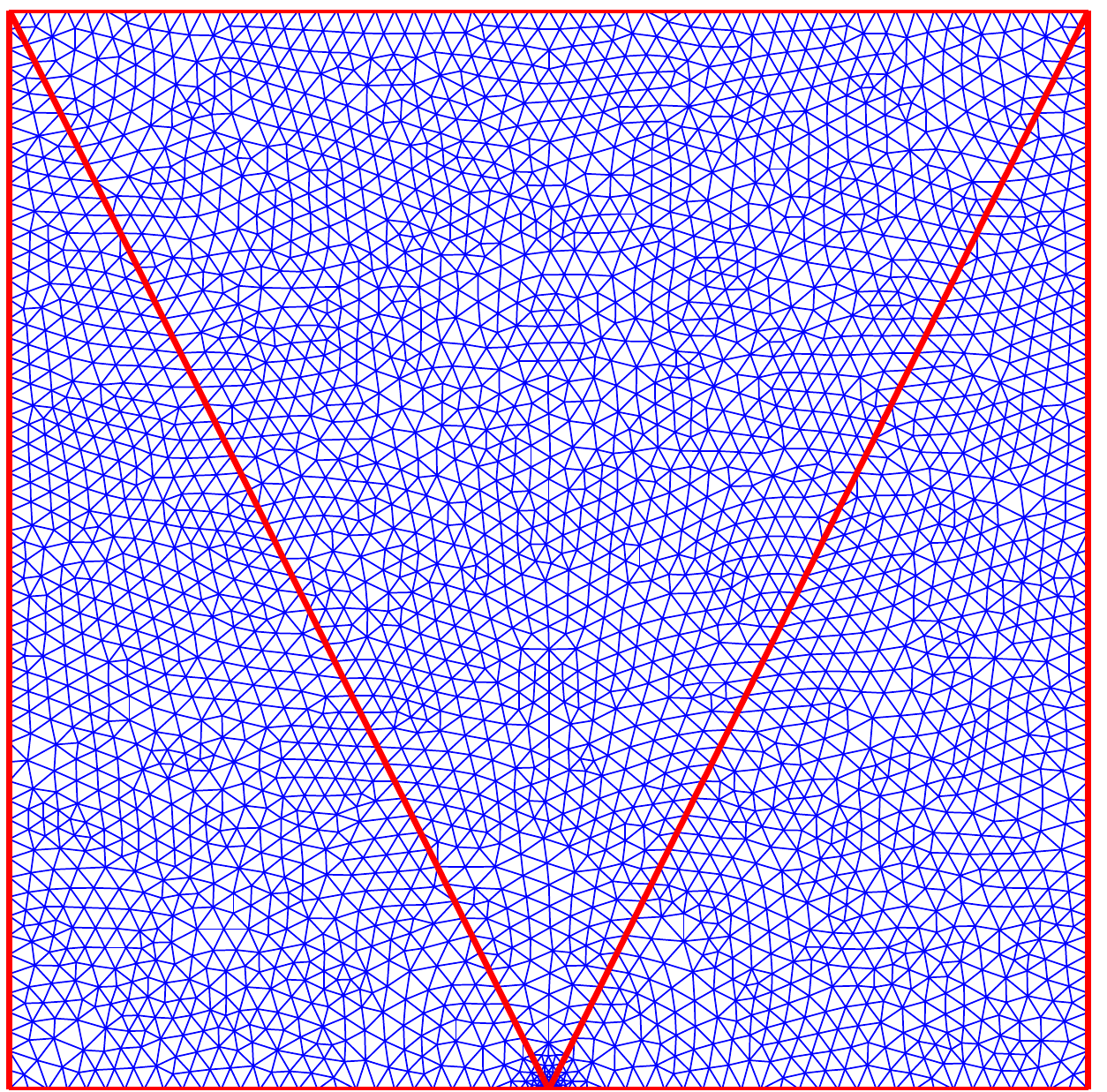}
\caption{The center crack problem: Domain composition and FE mesh}
\label{fig:ex2_mesh}
\end{figure}

We then consider a FE approximation with a mesh contains $n_{\rm node} = 3257$ nodes and $n_{\rm elem} = 6276$ $P_1$ elements, which corresponds to $\calN = 6422$ degrees of freedoms; the mesh is refined around the crack tip in order to give a good approximation for the (singular) solution near this region as shown in Figure~\ref{fig:ex2_mesh}.

We present in Table~\ref{tab:ex2_tab} the convergence results for the ``compliant'' output $s(\bfmu)$: the RB error bounds and effectivities as a function of $N$. The error bound reported, $\calE_N = \Delta^s_N(\bfmu)/|s_N(\bfmu)|$ is the maximum of the relative error bound over a random test sample $\Xi_{\rm test}$ of size $n_{\rm test} = 200$. We denote by $\overline{\eta}_N^s$ the average of the effectivity $\eta_N^s(\bfmu)$ over $\Xi_{\rm test}$. We observe that the effectivity average is very sharp, and of order $O(10)$.
\begin{table}
\centering
\begin{tabular}{|c||c|c|}
    \hline
    $N$ & $\calE_N$ & $\overline{\eta}_N^s$ \\
    \hline
     5  &  2.73\texttt{E}-02 &   6.16 \\
    10  &  9.48\texttt{E}-04 &   8.47 \\
    20  &  5.71\texttt{E}-06 &   7.39 \\
    30  &  5.59\texttt{E}-08 &   7.01 \\
    40  &  8.91\texttt{E}-10 &   7.54 \\
    50  &  6.26\texttt{E}-11 &   8.32 \\
    \hline
\end{tabular}
\caption{The center crack problem: RB convergence}
\label{tab:ex2_tab}
\end{table}

We next define the ERR RB approximation $\widehat{G}_N(\bfmu)$ to our ``truth'' (FE) $\widehat{G}^\calN_{\rm FE}(\bfmu)$ and its associated ERR RB error $\Delta^{\widehat{G}}_N(\bfmu)$ by
\begin{eqnarray}\label{eqn:RB_SIF_err}
    \widehat{G}_N(\bfmu) &=& \frac{s_N(\bfmu) - \Delta_N^s(\bfmu+\delta\mu_1)}{\delta\mu_1}, \nonumber \\
    \Delta^{\widehat{G}}_N(\bfmu) &=& \frac{\Delta_N^s(\bfmu + \delta\mu_1) + \Delta_N^s(\bfmu)}{\delta\mu_1}.
\end{eqnarray}
It can be readily proven \cite{rozza08:ARCME} that our SIF RB error is a rigorous bound for the ERR RB prediction $\widehat{G}_N(\bfmu)$: $|\widehat{G}_N(\bfmu) - \widehat{G}^\calN_{\rm FE}(\bfmu)| \leq \Delta^{\widehat{G}}_N(\bfmu)$. It is note that the choice of $\delta\mu_1$ is not arbitrary: $\delta\mu_1$ needed to be small enough to provide a good FD approximation, while still provide a good ERR RB error bound \refeq{eqn:RB_SIF_err}. Here we choose $\delta\mu_1 = 1\texttt{E}-03$.

We then can define the SIF RB approximation $\widehat{\rm SIF}_N(\bfmu)$ to our ``truth'' (FE) $\widehat{\rm SIF}^\calN_{\rm FE}(\bfmu)$ and its associated SIF RB error estimation $\Delta^{\widehat{\rm SIF}}_N(\bfmu)$ as
\begin{eqnarray*}
\widehat{\rm SIF}_N(\bfmu) &=& \frac{1}{2\sqrt{1-\nu^2}}\bigg\{\sqrt{\widehat{G}_N(\bfmu) + \Delta^{\widehat{G}}_N(\bfmu)}+\sqrt{\widehat{G}_N(\bfmu) - \Delta^{\widehat{G}}_N(\bfmu)}\bigg\}, \\
\Delta^{\widehat{\rm SIF}}_N(\bfmu) &=& \frac{1}{2\sqrt{1-\nu^2}}\bigg\{\sqrt{\widehat{G}_N(\bfmu) + \Delta^{\widehat{G}}_N(\bfmu)}-\sqrt{\widehat{G}_N(\bfmu) - \Delta^{\widehat{G}}_N(\bfmu)}\bigg\}.
\end{eqnarray*}
It is readily proven in \cite{huynh07:ijnme} that $|\widehat{\rm SIF}_N(\bfmu) - \widehat{\rm SIF}^\calN_{\rm FE}(\bfmu)| \leq \Delta^{\widehat{\rm SIF}}_N(\bfmu)$.

\begin{figure} [htbp]
\centering
\includegraphics[scale=0.4]{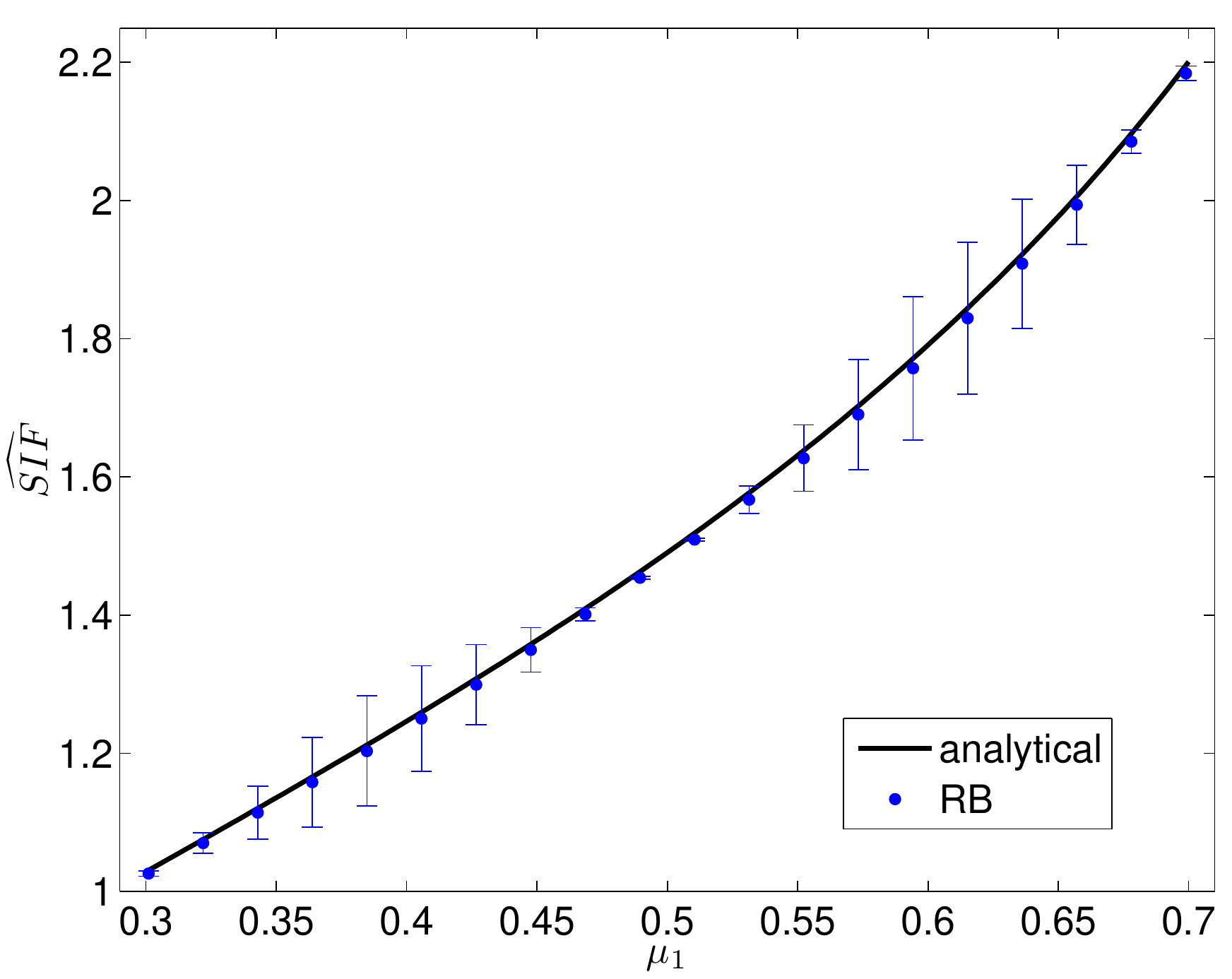}
\caption{The center crack problem: SIF solution for $N=15$}
\label{fig:ex2_SIF15}
\end{figure}

\begin{figure} [htbp]
\centering
\includegraphics[scale=0.4]{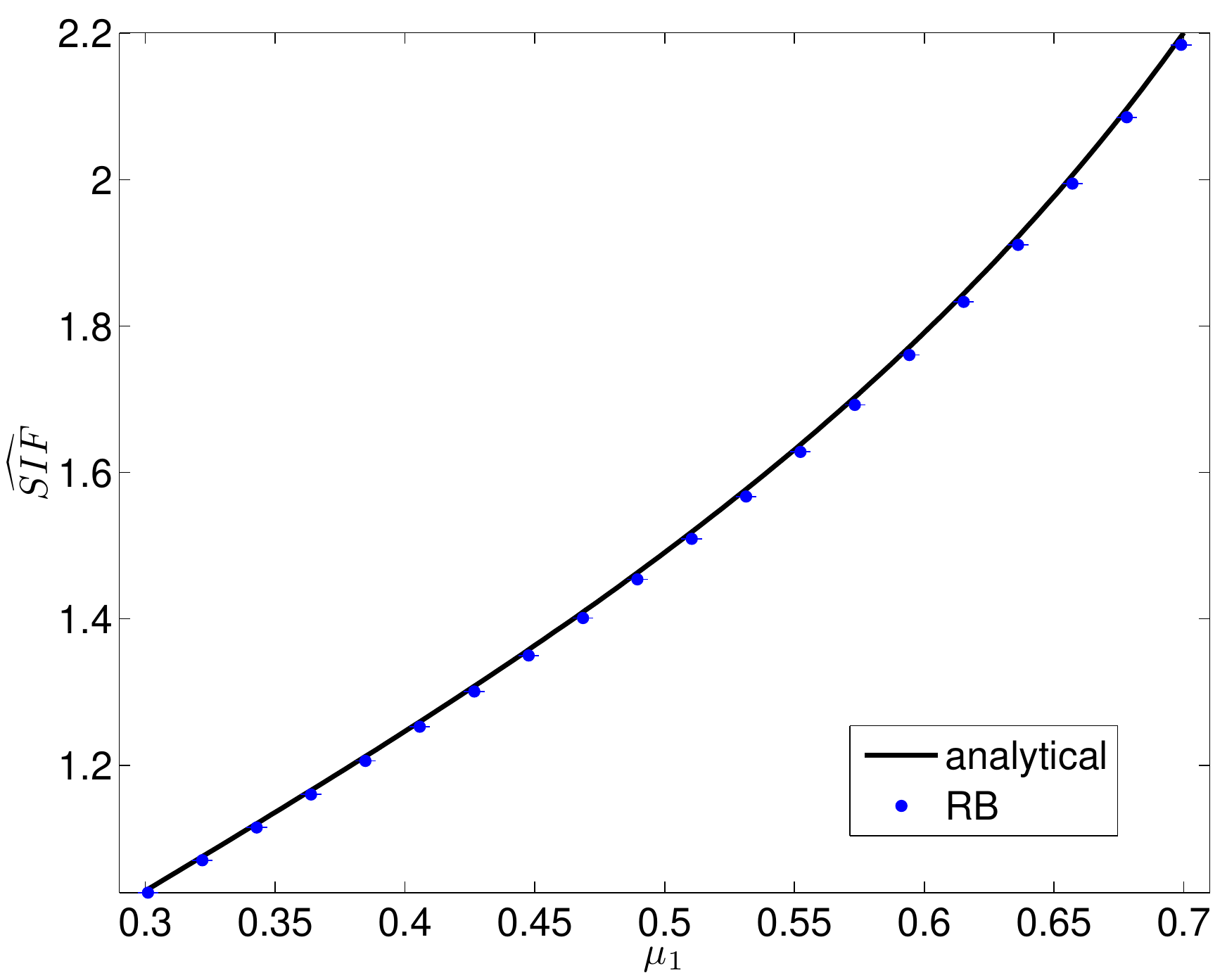}
\caption{The center crack problem: SIF solution for $N=30$}
\label{fig:ex2_SIF30}
\end{figure}

We plot the SIF RB results $\widehat{\rm SIF}(\bfmu)$ with error bars correspond to $\Delta^{\widehat{\rm SIF}}_N(\bfmu)$, and the analytical results $\widehat{\rm SIF}(\bfmu)$ \cite{murakami01:SIFhandbook}
in Figure~\ref{fig:ex2_SIF15} for the case $\mu_1 \in [0.3,0.7]$, $\mu_2 = 2.0$ for $N = 15$. It is observed that the RB error is large since the small number of basis $N = 15$ does not compromise the small $\delta\mu_1 = 1\texttt{E}-03$ value.
We next plot, in Figure~\ref{fig:ex2_SIF30}, SIF RB results and error for the same $\bfmu$ range as in Figure~\ref{fig:ex2_SIF15}, but for $N = 30$. It is observed now that the SIF RB error is significantly improved -- thanks to the better RB approximation that compensates the small value $\delta\mu_1$. We also want to point out that, in both Figure~\ref{fig:ex2_SIF15} and Figure~\ref{fig:ex2_SIF30}, it is clearly shown that our RB SIF error is not a \emph{rigorous} bound for the \emph{exact} SIF values $\widehat{\rm SIF}(\bfmu)$ but rather is a \emph{rigorous} bound for the ``truth'' (FE) approximation $\widehat{\rm SIF}^\calN_{\rm FE}(\bfmu)$. It is shown, however, that FE SIF approximation (which is considered in Figure~\ref{fig:ex2_SIF30} thanks to the negligible RB error) are of good quality compared with the exact SIF. The VCE in this case works quite well, however it is not suitable for complicate crack settings. In such cases, other SIF calculation methods and appropriate RB approximations might be preferable \cite{huynh07:ijnme, huynh07:_phd_thesis}.

As regards computational times, a RB online evaluation $\bfmu \rightarrow (\widehat{\rm SIF}_N(\bfmu),\Delta_N^{\widehat{\rm SIF}}(\bfmu)$ requires just $t_{\rm RB} (= 25 \times) = 50$(ms) for $N = 40$; while the FE solution $\bfmu \rightarrow \widehat{\rm SIF}^\calN_{\rm FE}(\bfmu)$ requires $t_{\rm FE} (= 7 \times 2) = 14$(s): thus our RB online evaluation takes only $0.36\%$ of the FEM computational cost.

\subsection{The composite unit cell problem}

\begin{figure} [htbp]
\centering
\includegraphics[scale=0.5]{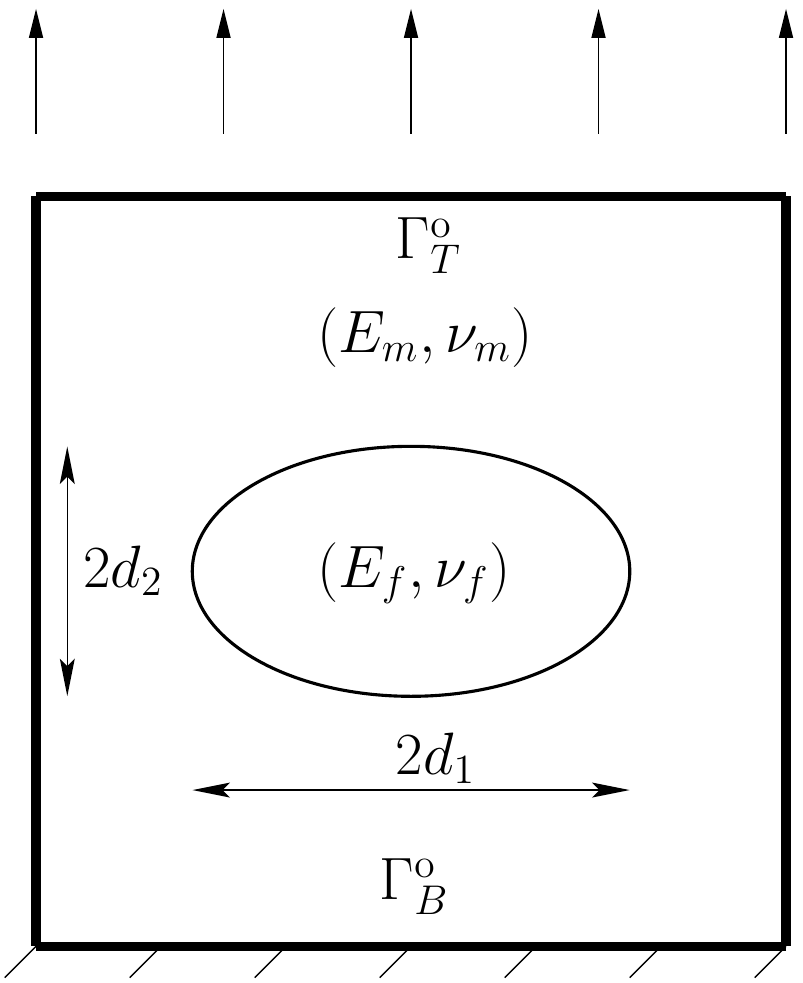}
\caption{The composite unit cell problem}
\label{fig:ex3_model}
\end{figure}

We consider a unit cell contains an ellipse region as shown in Figure~\ref{fig:ex3_model}. We apply (clamped) Dirichlet boundary conditions on the bottom of the cell $\Gamma^{\rm o}_B$ and (unit tension) non-homogeneous Neumann boundary conditions on $\Gamma^{\rm o}_T$. We denote the two semimajor axis and semiminor axis of the ellipse region as $d_1$ and $d_2$, respectively. We assume plane stress isotropic materials: the material properties of the matrix (outside of the region) is given by $(E_m,\nu_m) = (1,0.3)$, and the material properties of the ellipse region is given by $(E_f,\nu_f) = (E_f,0.3)$. Our output of interest is the integral of normal displacement ($u_1$) over $\Gamma^{\rm o}_T$. We note our output of interest is thus ``compliant''.

We consider $P=3$ parameters $\bfmu = [\mu_1,\mu_2,\mu_3] \equiv [d_1,d_2,E_f]$. The parameter domain is chosen as $\calD = [0.8, 1.2] \times [0.8,1.2] \times [0.2,5]$. Note that the third parameter (the Young modulus of the ellipse region) can represent the ellipse region from an ``inclusion'' (with softer Young's modulus $E_f<E_m (= 1)$) to a ``fiber'' (with stiffer Young's modulus $E_f>E_m (= 1)$).

\begin{figure} [htbp]
\centering
\includegraphics[scale=0.3]{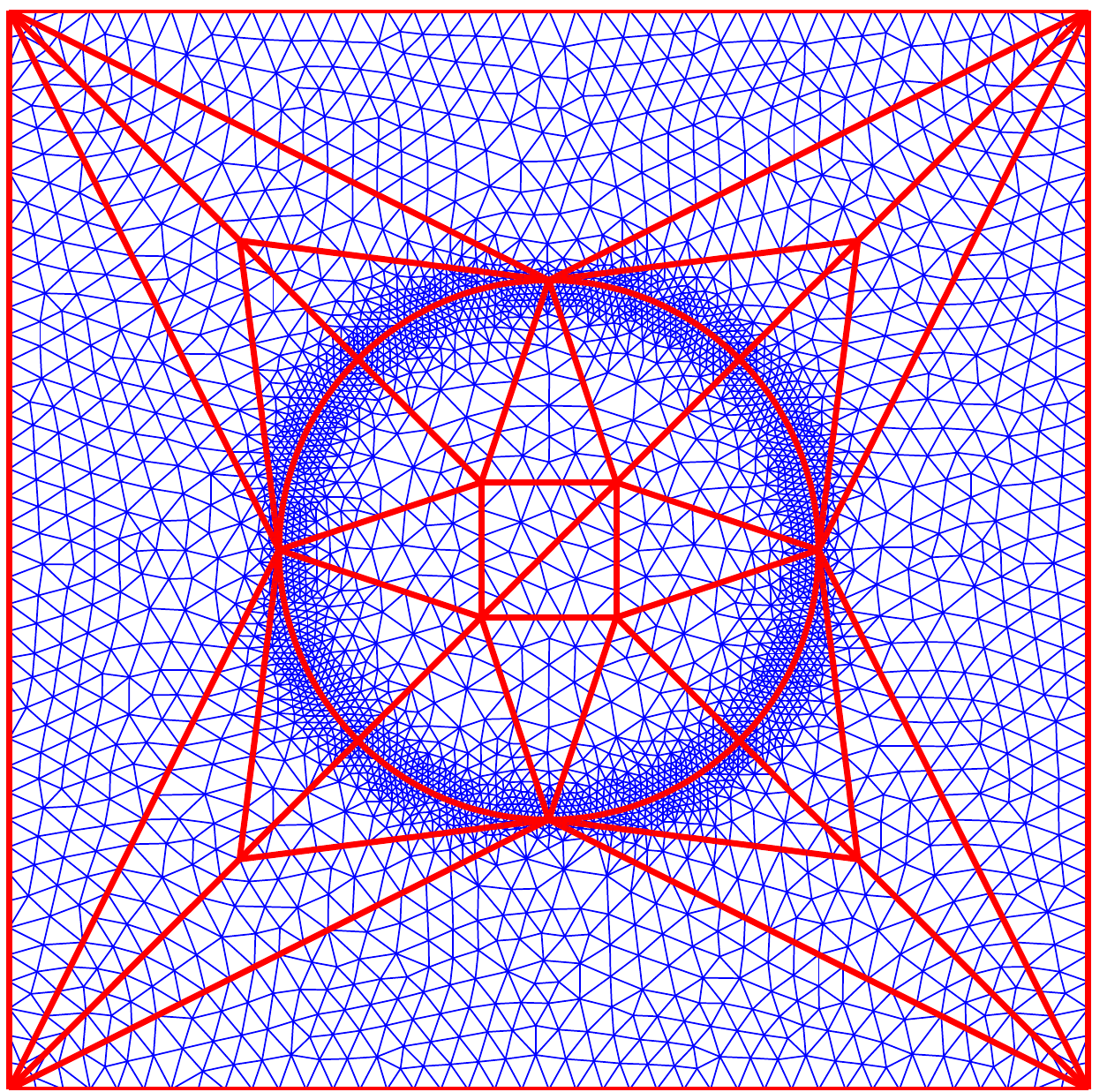}
\caption{The composite unite cell problem: Domain composition and FE mesh}
\label{fig:ex3_mesh}
\end{figure}

We then choose $\bfmu_{\rm ref} = [1.0,1.0,1.0]$ and apply the domain decomposition \cite{rozza08:ARCME} and obtain $L_{\rm reg} = 34$ subdomains, in which $16$ subdomains are the general ``curvy triangles'' ($8$ inward ``curvy triangles'' and $8$ outward curvy ``triangles'') as shown in Figure~\ref{fig:ex3_mesh}. However, despite the large number of ``curvy triangles'' in the domain decomposition, it is observed that almost all transformations are congruent, hence we expected a small number of $Q^a$ than (says), that of the ``arc-cantilever beam'' example, in which all the subdomains transformations are different. Indeed, we recover our affine forms with $Q^a = 30$ and $Q^f = 1$, note that $Q^a$ is relatively small for such a complex domain decomposition thanks to our efficient symbolic manipulation ``collapsing'' technique and those congruent ``curvy triangles''.

We next consider a FE approximation where the mesh contains $n_{\rm node} = 3906$ nodes and $n_{\rm elem} = 7650$ $P_1$ elements, which corresponds to $\calN = 7730$ degrees of freedoms. The mesh is refined around the interface of the matrix and the inclusion/fiber.

We then apply the RB approximation.  We present in Table~\ref{tab:ex3_tab} our convergence results: the RB error bounds and effectivities as a function of $N$. The error bound reported, $\calE_N = \Delta^s_N(\bfmu)/|s_N(\bfmu)|$ is the maximum of the relative error bound over a random test sample $\Xi_{\rm test}$ of size $n_{\rm test} = 200$. We denote by $\overline{\eta}_N^s$ the average of the effectivity $\eta_N^s(\bfmu)$ over $\Xi_{\rm test}$. We observe that our effectivity average is of order $O(10)$.
\begin{table}
\centering
\begin{tabular}{|c||c|c|}
    \hline
    $N$ & $\calE_N$ & $\overline{\eta}_N^s$ \\
    \hline
     5 &  9.38\texttt{E}-03 &   8.86 \\
    10 &  2.54\texttt{E}-04 &   7.18 \\
    15 &  1.37\texttt{E}-05 &   5.11 \\
    20 &  3.91\texttt{E}-06 &   9.74 \\
    25 &  9.09\texttt{E}-07 &   6.05 \\
    30 &  2.73\texttt{E}-07 &   10.64 \\
    35 &  9.00\texttt{E}-08 &   10.17 \\
    40 &  2.66\texttt{E}-08 &   10.35 \\
    \hline
\end{tabular}
\caption{The composite unit cell problem: RB convergence}
\label{tab:ex3_tab}
\end{table}

As regards computational times, a RB online evaluation $\bfmu \rightarrow (s_N(\bfmu),\Delta_N^s(\bfmu))$ requires just $t_{\rm RB} = 66$(ms) for $N = 30$; while the FE solution $\bfmu \rightarrow s^\calN(\bfmu)$ requires approximately $t_{\rm FE} = 8$(s): thus our RB online evaluation is just $0.83\%$ of the FEM computational cost.

\subsection{The multi-material plate problem}

\begin{figure} [htbp]
\centering
\includegraphics[scale=0.5]{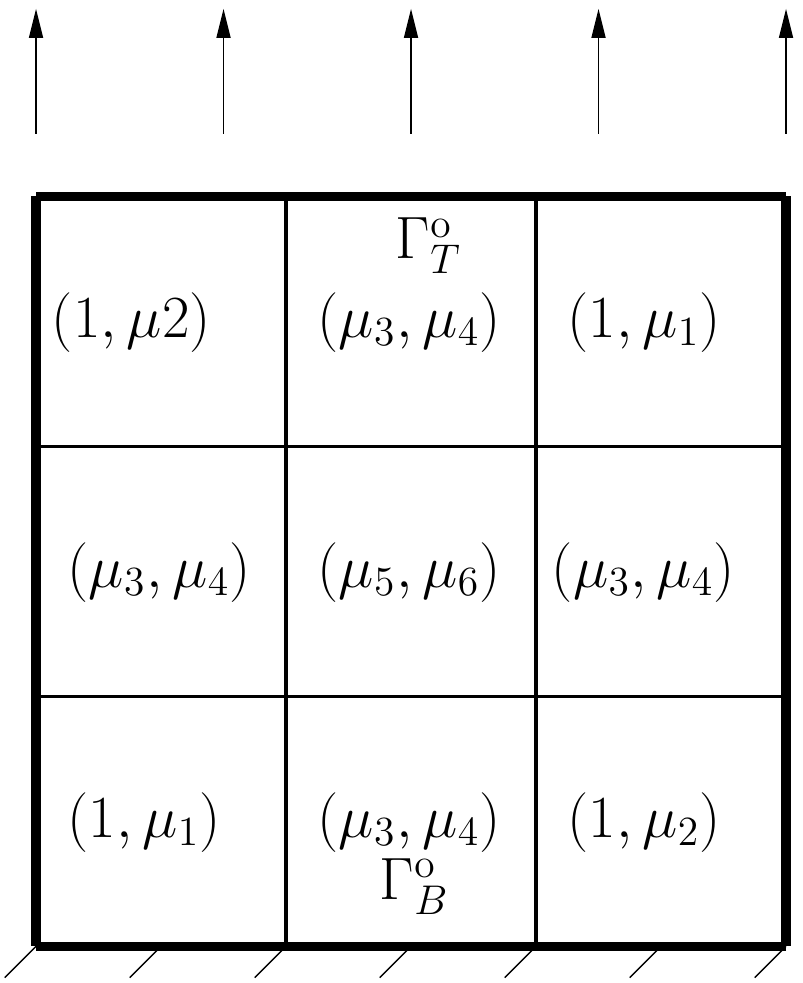}
\caption{The multi-material problem}
\label{fig:ex4_model}
\end{figure}

We consider a unit cell divided into $9$ square subdomains of equal size as shown in Figure~\ref{fig:ex4_model}. We apply (clamped) Dirichlet boundary conditions on the bottom of the cell $\Gamma^{\rm o}_B$ and (unit tension) non-homogeneous Neumann boundary conditions on $\Gamma^{\rm o}_T$. We consider orthotropic plane stress materials: the Young's modulus properties for all $9$ subdomains are given in Figure, the Poisson's ratio is chosen as $\nu_{12,i} =0.3$, $i = 1,\ldots,9$ and $\nu_{21,i}$ is determined by \refeq{eqn:ortho_plane_stress}. The shear modulus is chosen as a function of the two Young's moduli as in \refeq{eqn:ortho_shear_modulus} for all $9$ subdomains. All material axes are aligned with the coordinate system (and loading). Our output of interest is the integral of normal displacement ($u_1$) over $\Gamma^{\rm o}_T$, which represents the average normal displacement on $\Gamma^{\rm o}_T$. We note our output of interest is thus ``compliant''.

We consider $P=6$ parameters $\bfmu = [\mu_1,\ldots,\mu_6]$, correspond to the six Young's moduli values as shown in Figure~\ref{fig:ex4_model} (the two Young's moduli for each subdomain are shown in those brackets). The parameter domain is chosen as $\calD = [0.5, 2.0]^6$.

We then apply the domain decomposition \cite{rozza08:ARCME} and obtain $L_{\rm reg} = 18$ subdomains. Despite the large $L_{\rm reg}$ number of domains, there is no geometric transformation in this case. We recover our affine forms with $Q^a = 12$, $Q^f = 1$, note that all $Q^a$ are contributed from all the Young's moduli since there is no geometric transformation involved. Moreover, it is observed that the bilinear form can be, in fact, classified as a ``parametrically coercive'' one \cite{patera07:book}.

\begin{figure} [htbp]
\centering
\includegraphics[scale=0.3]{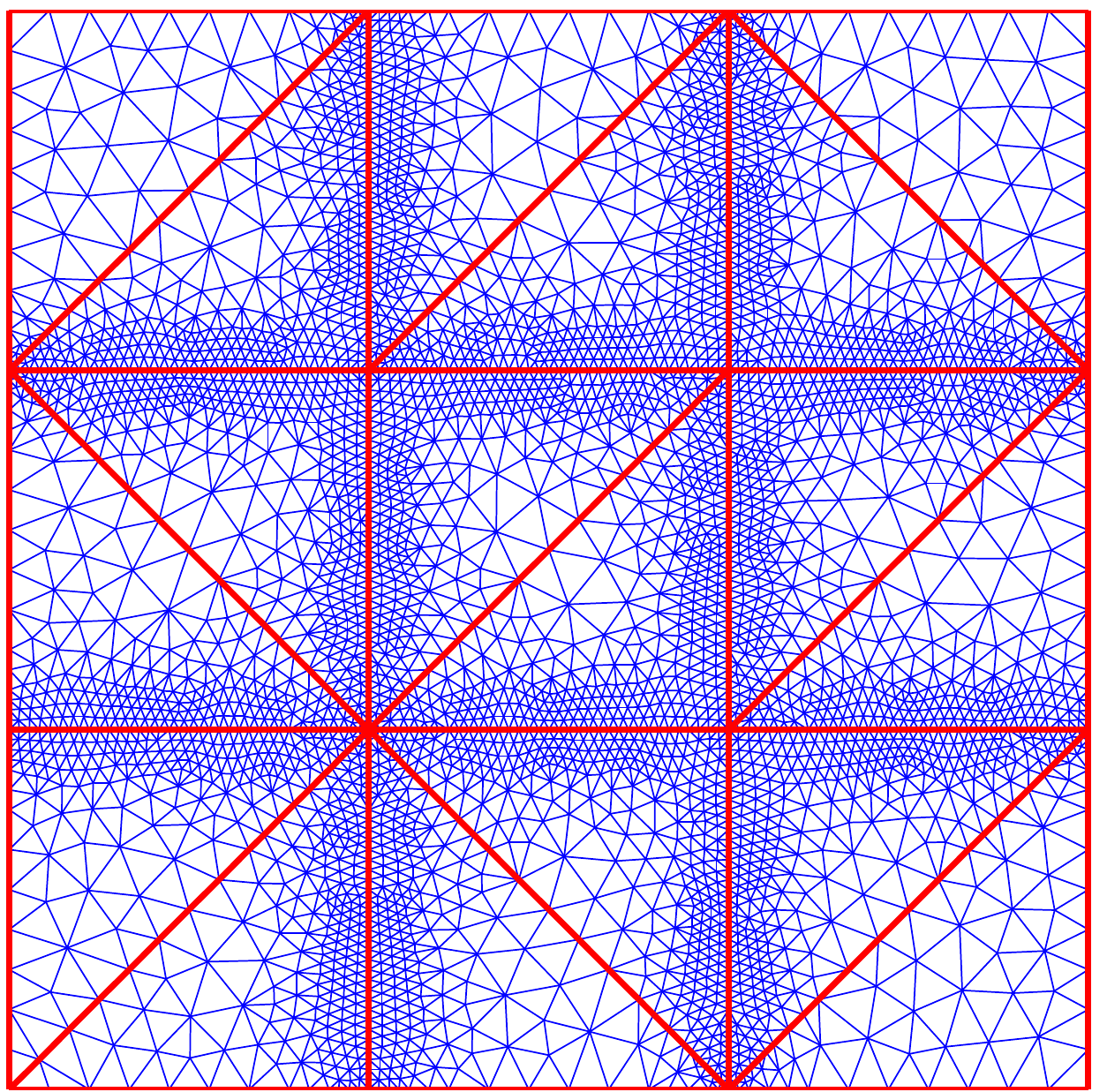}
\caption{The multi-material problem: Domain composition and FE mesh}
\label{fig:ex4_mesh}
\end{figure}

We next consider a FE approximation where the mesh contains $n_{\rm node} = 4098$ nodes and $n_{\rm elem} = 8032$ $P_1$ elements, which corresponds to $\calN = 8112$ degrees of freedoms. The mesh is refined around all the interfaces between different subdomains as shown in Figure~\ref{fig:ex4_mesh}.

We then apply the RB approximation.  We present in Table~\ref{tab:ex4_tab} our convergence results: the RB error bounds and effectivities as a function of $N$. The error bound reported, $\calE_N = \Delta^s_N(\bfmu)/|s_N(\bfmu)|$ is the maximum of the relative error bound over a random test sample $\Xi_{\rm test}$ of size $n_{\rm test} = 200$. We denote by $\overline{\eta}_N^s$ the average of the effectivity $\eta_N^s(\bfmu)$ over $\Xi_{\rm test}$. We observe that our effectivity average is of order $O(10)$.
\begin{table}
\centering
\begin{tabular}{|c||c|c|}
    \hline
    $N$ & $\calE_N$ & $\overline{\eta}_N^s$ \\
    \hline
     5 &  1.01\texttt{E}-02 &   8.11 \\
    10 &  1.45\texttt{E}-03 &   11.16 \\
    20 &  3.30\texttt{E}-04 &   11.47 \\
    30 &  1.12\texttt{E}-04 &   12.59 \\
    40 &  2.34\texttt{E}-05 &   11.33 \\
    50 &  9.85\texttt{E}-06 &   12.90 \\
    \hline
\end{tabular}
\caption{The multi-material problem: RB convergence}
\label{tab:ex4_tab}
\end{table}

As regards computational times, a RB online evaluation $\bfmu \rightarrow (s_N(\bfmu),\Delta_N^s(\bfmu))$ requires just $t_{\rm RB} = 33$(ms) for $N = 40$; while the FE solution $\bfmu \rightarrow s^\calN(\bfmu)$ requires $t_{\rm FE} = 8.1$(s): thus the RB online evaluation is just $0.41\%$ of the FEM computational cost.

\subsection{The woven composite beam problem}

\begin{figure} [htbp]
\centering
\includegraphics[scale=0.5]{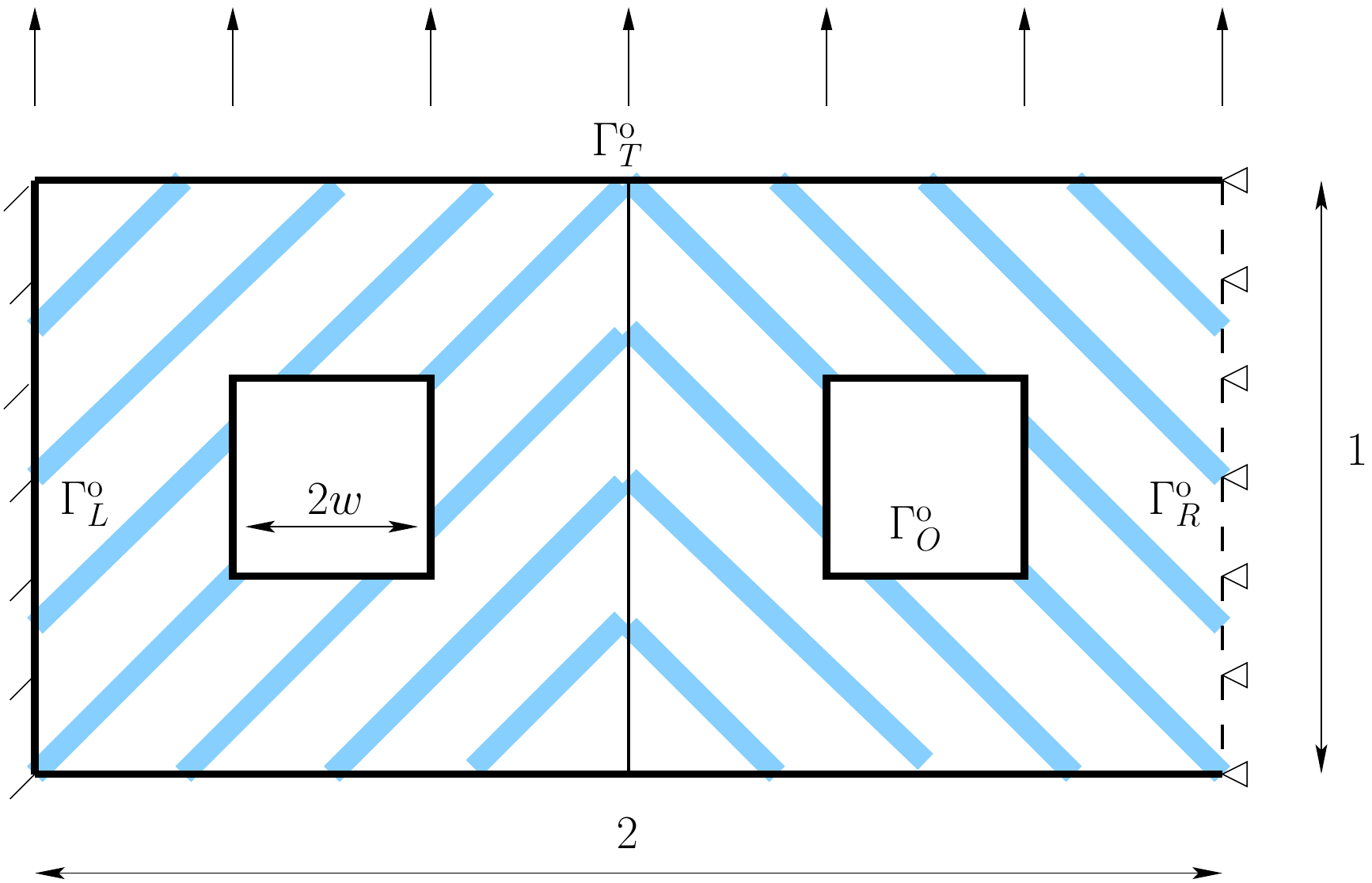}
\caption{The woven composite beam problem}
\label{fig:ex5_model}
\end{figure}

We consider a composite cantilever beam as shown in Figure~\ref{fig:ex5_model}. The beam is divided into two regions, each with a square hole in the center of (equal) size $2w$. We apply (clamped) Dirichlet boundary conditions on the left side of the beam $\Gamma^{\rm o}_L$, (symmetric about the $x^{\rm o}_1$ direction) Dirichlet boundary conditions on the right side of the beam $\Gamma^{\rm o}_R$, and (unit tension) non-homogeneous Neumann boundary conditions on the top side $\Gamma^{\rm o}_T$. We consider the same orthotropic plane stress materials for both regions: $(E_1, E_2) = (1, E_2)$, $\nu_{12} = 0.3$, $\nu_{21}$ is determined by \refeq{eqn:ortho_plane_stress} and the shear modulus $G_{12}$ is given by \refeq{eqn:ortho_shear_modulus}. The material axes of both regions are not aligned with the coordinate system and loading: the angles of the the material axes and the coordinate system of the first and second region are $\theta$ and $-\theta$, respectively. The setting represents a ``woven'' composite material across the beam horizontally. Our output of interest is the integral of the normal displacement ($u_1$) over the boundary $\Gamma^{\rm o}_O$. We note our output of interest is thus ``non-compliant''.

We consider $P=3$ parameters $\bfmu = [\mu_1,\mu_2,\mu_3] \equiv [w, E_2, \theta]$. The parameter domain is chosen as $\calD = [1/6, 1/12] \times [1/2, 2] \times [-\pi/4, \pi/4]$.

We then apply the domain decomposition \cite{rozza08:ARCME} and obtain $L_{\rm reg} = 32$ subdomains, note that all subdomains transformations are just simply translations due to the ``added control points'' along the external (and interface) boundaries strategy \cite{rozza08:ARCME}. We recover the affine forms with $Q^a = 19$, $Q^f = 2$, and $Q^\ell = 1$.

\begin{figure} [htbp]
\centering
\includegraphics[scale=0.5]{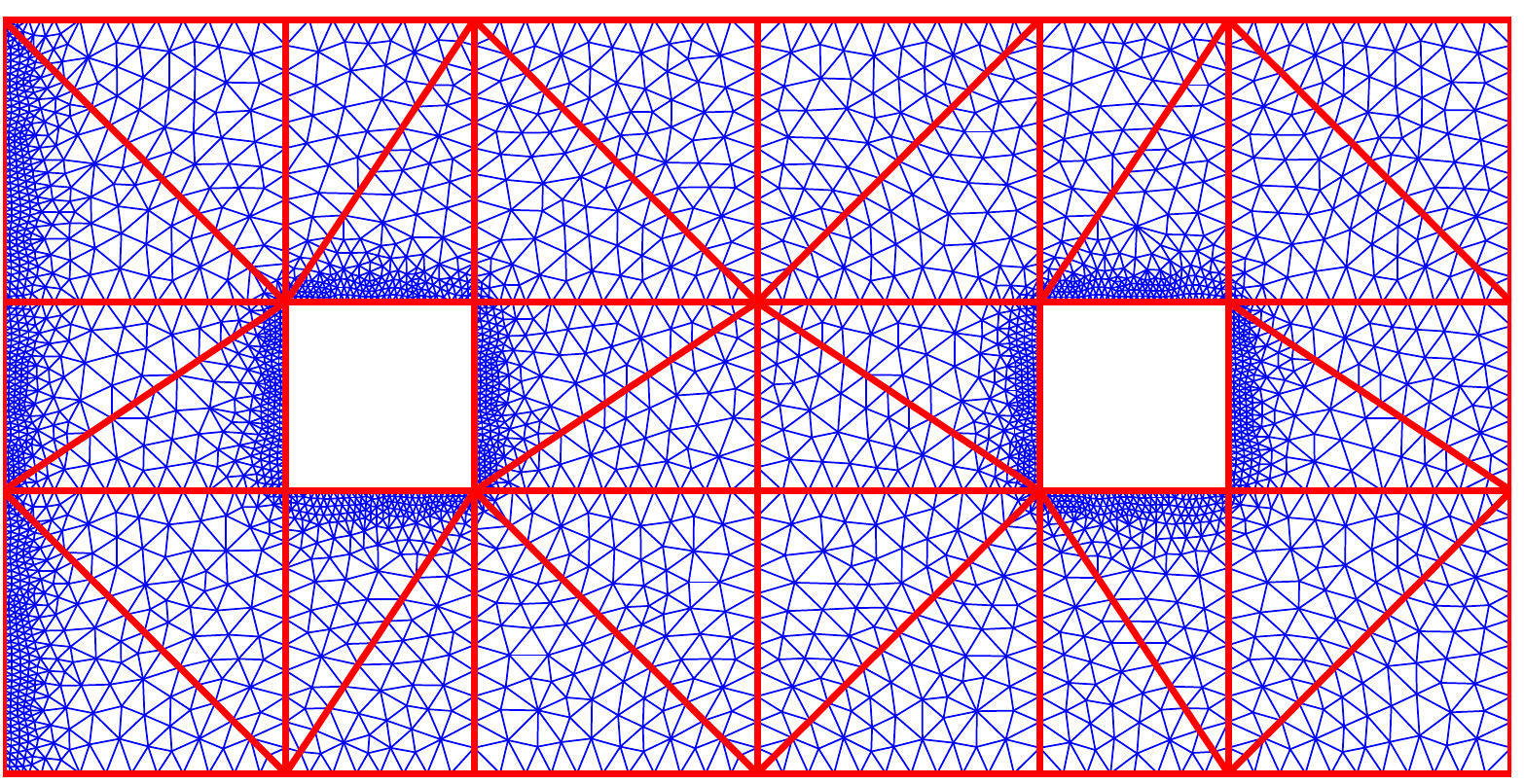}
\caption{The woven composite beam problem: Domain composition and FE mesh}
\label{fig:ex5_mesh}
\end{figure}

We next consider a FE approximation where the mesh contains $n_{\rm node} = 3569$ nodes and $n_{\rm elem} = 6607$ $P_1$ elements, which corresponds to $\calN = 6865$ degrees of freedoms. The mesh is refined around the holes, the interfaces between the two regions, and the clamped boundary as shown in Figure~\ref{fig:ex5_mesh}.

We then apply the RB approximation. We present in Table~\ref{tab:ex5_tab} our convergence results: the RB error bounds and effectivities as a function of $N$. The error bound reported, $\calE_N = \Delta^s_N(\bfmu)/|s_N(\bfmu)|$ is the maximum of the relative error bound over a random test sample $\Xi_{\rm test}$ of size $n_{\rm test} = 200$. We denote by $\overline{\eta}_N^s$ the average of the effectivity $\eta_N^s(\bfmu)$ over $\Xi_{\rm test}$. We observe that our effectivity average is of order $O(5-25)$.
\begin{table}
\centering
\begin{tabular}{|c||c|c|}
    \hline
    $N$ & $\calE_N$ & $\overline{\eta}_N^s$ \\
    \hline
     4 &  4.64\texttt{E}-02 &   22.66 \\
     8 &  1.47\texttt{E}-03 &    7.39 \\
    12 &  2.35\texttt{E}-04 &    9.44 \\
    16 &  6.69\texttt{E}-05 &   14.29 \\
    20 &  1.31\texttt{E}-05 &   11.41 \\
    \hline
\end{tabular}
\caption{The woven composite beam problem: RB convergence}
\label{tab:ex5_tab}
\end{table}

As regards computational times, a RB online evaluation $\bfmu \rightarrow (s_N(\bfmu),\Delta_N^s(\bfmu))$ requires just $t_{\rm RB} = 40$(ms) for $N = 20$; while the FE solution $\bfmu \rightarrow s^\calN(\bfmu)$ requires $t_{\rm FE} = 7.5$(s): thus our RB online evaluation is just $0.53\%$ of the FEM computational cost.

\subsection{The closed vessel problem}

\begin{figure} [htbp]
\centering
\includegraphics[scale=0.3]{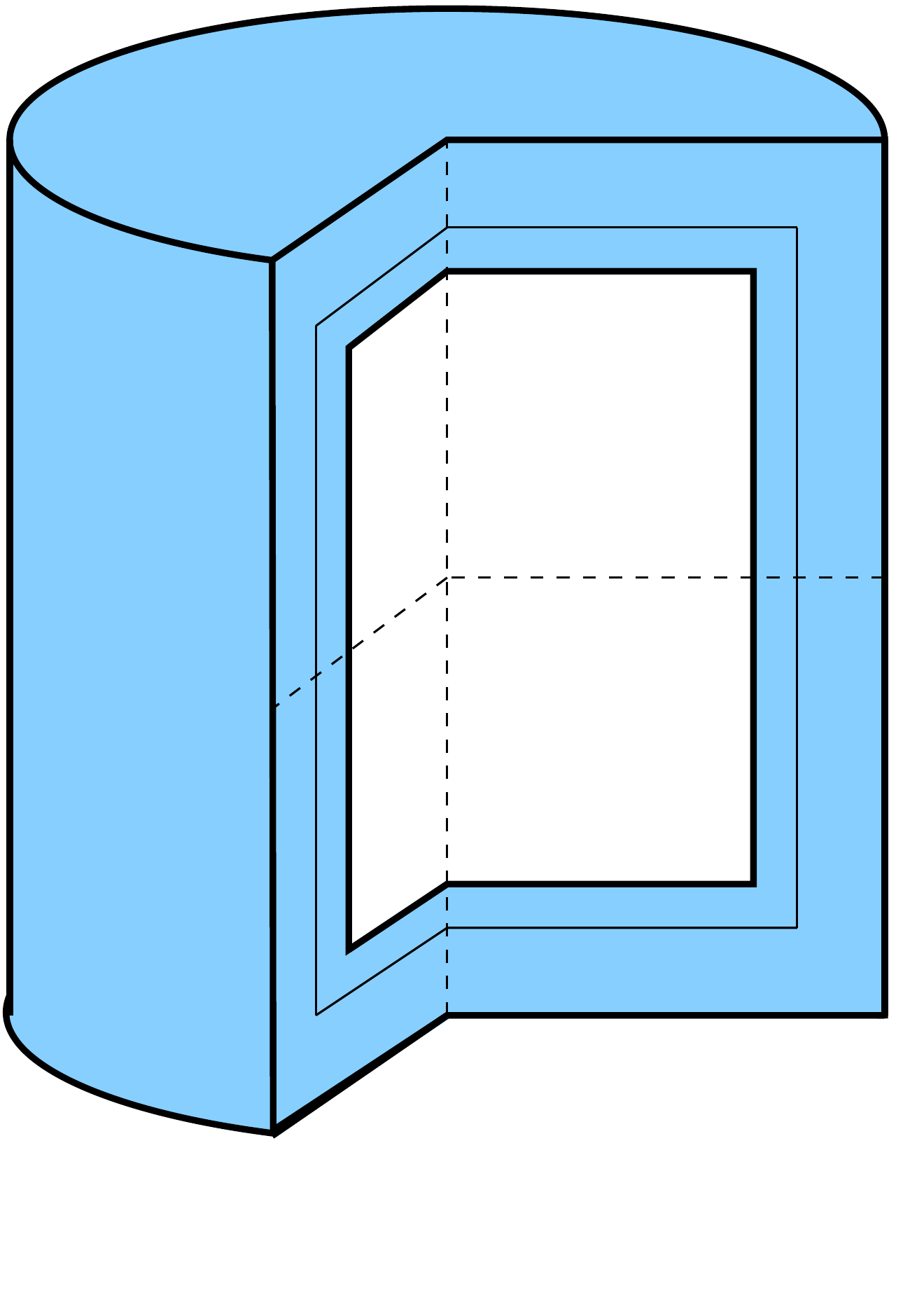}
\caption{The closed vessel problem}
\label{fig:ex6_modelfull}
\end{figure}

\begin{figure} [htbp]
\centering
\includegraphics[scale=0.5]{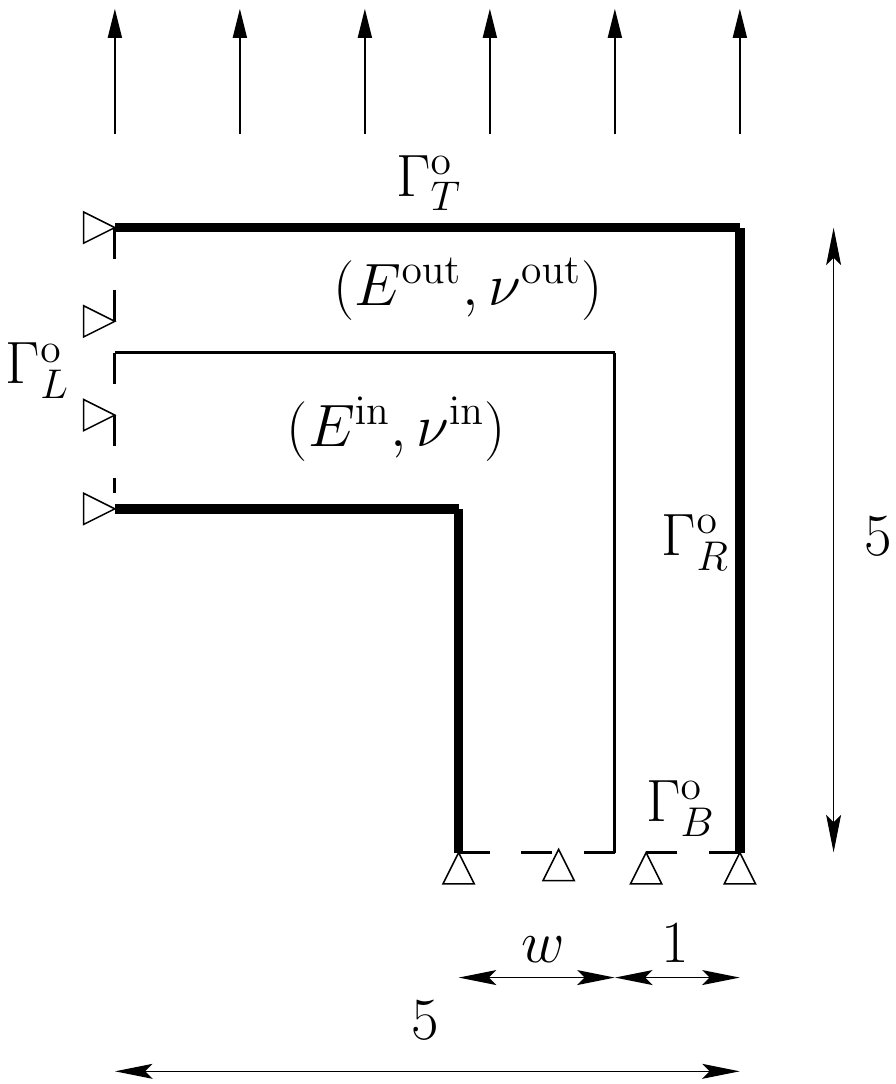}
\caption{The closed vessel problem}
\label{fig:ex6_model}
\end{figure}

We consider a closed vessel under tension at both ends as shown in Figure~\ref{fig:ex6_modelfull}.
The vessel is axial symmetric about the $x^{\rm o}_2$ axis, and symmetric about the $x^{\rm o}_1$ axis, hence we only consider a representation ``slice'' by our axisymmetric formulation as shown in Figure~\ref{fig:ex6_model}. The vessel is consists of two layered, the outer layer is of fixed width $w^{\rm out} = 1$, while the inner layer is of width $w^{\rm in} = w$. The material properties of the inner layer and outer layer are given by $(E^{\rm in},\nu) = (E^{\rm in},0.3)$ and $(E^{\rm out},\nu) = (1,0.3)$, respectively. We apply (symmetric about the $x^{\rm o}_2$ direction) Dirichlet boundary conditions on the bottom boundary of the model $\Gamma^{\rm o}_B$, (symmetric about the $x^{\rm o}_1$ direction) Dirichlet boundary conditions on the left boundary of the model $\Gamma^{\rm o}_L$ and (unit tension) non-homogeneous Neumann boundary conditions on the top boudanry  $\Gamma^{\rm o}_T$. Our output of interest is the integral of the axial displacement ($u_r$) over the right boundary $\Gamma^{\rm o}_R$. We note our output of interest is thus ``non-compliant''.

We consider $P=2$ parameters $\bfmu = [\mu_1,\mu_2] \equiv [w, E^{\rm in}]$. The parameter domain is chosen as $\calD = [0.1, 1.9] \times [0.1, 10]$.

\begin{figure} [htbp]
\centering
\includegraphics[scale=0.3]{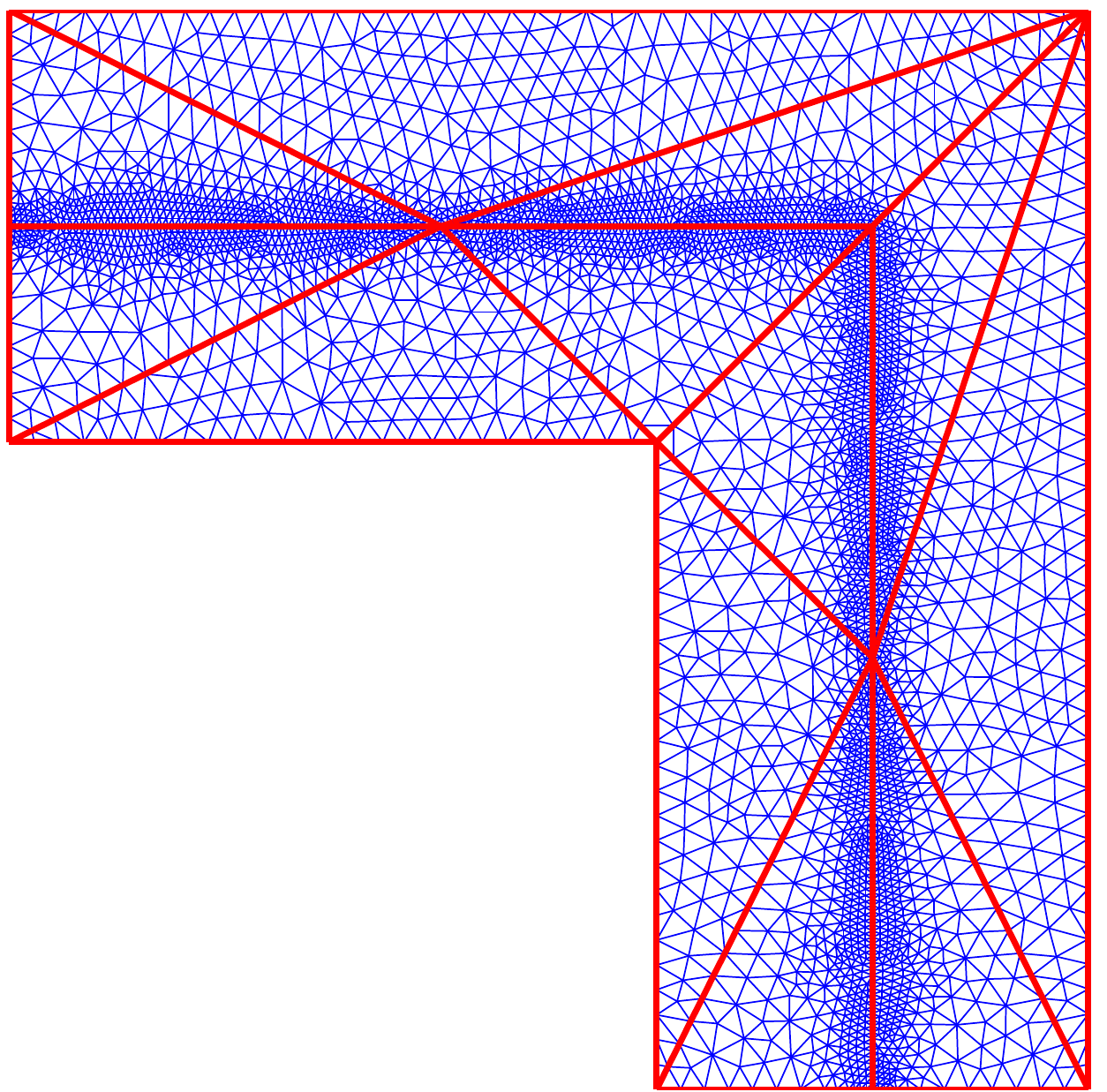}
\caption{The closed vessel problem: Domain composition and FE mesh}
\label{fig:ex6_mesh}
\end{figure}

We then apply the domain decomposition \cite{rozza08:ARCME} and obtain $L_{\rm reg} = 12$ subdomains as shown in Figure~\ref{fig:ex6_mesh}. We recover our affine forms with $Q^a = 47$, $Q^f = 1$, and $Q^\ell = 1$. Despite the small number of parameter (and seemingly simple transformations), $Q^a$ is large in this case. A major contribution to $Q^a$ come from the expansion of the terms $x^{\rm o}_1$ in the effective elastic tensor $[\bfS]$, which appeared due to the geometric transformation of the inner layer.

We next consider a FE approximation where the mesh contains $n_{\rm node} = 3737$ nodes and $n_{\rm elem} = 7285$ $P_1$ elements, which corresponds to $\calN = 7423$ degrees of freedoms. The mesh is refined around the interfaces between the two layers.

We then apply the RB approximation.  We present in Table~\ref{tab:ex6_tab} convergence results: the RB error bounds and effectivities as a function of $N$. The error bound reported, $\calE_N = \Delta^s_N(\bfmu)/|s_N(\bfmu)|$ is the maximum of the relative error bound over a random test sample $\Xi_{\rm test}$ of size $n_{\rm test} = 200$. We denote by $\overline{\eta}_N^s$ the average of the effectivity $\eta_N^s(\bfmu)$ over $\Xi_{\rm test}$. We observe that our effectivity average is of order $O(50-120)$, which is quite large, however it is not surprising since our output is ``non-compliant''.
\begin{table}
\centering
\begin{tabular}{|c||c|c|}
    \hline
    $N$ & $\calE_N$ & $\overline{\eta}_N^s$ \\
    \hline
    10 &  7.12\texttt{E}-02 &  56.01 \\
    20 &  1.20\texttt{E}-03 &  111.28 \\
    30 &  3.96\texttt{E}-05 &  49.62 \\
    40 &  2.55\texttt{E}-06 &  59.96 \\
    50 &  5.70\texttt{E}-07 &  113.86 \\
    60 &  5.90\texttt{E}-08 &  111.23 \\
    70 &  6.95\texttt{E}-09 &  77.12 \\
    \hline
\end{tabular}
\caption{The closed vessel problem: RB convergence}
\label{tab:ex6_tab}
\end{table}

As regards computational times, a RB online evaluation $\bfmu \rightarrow (s_N(\bfmu),\Delta_N^s(\bfmu))$ requires just $t_{\rm RB} = 167$(ms) for $N = 40$; while the FE solution $\bfmu \rightarrow s^\calN(\bfmu)$ requires $t_{\rm FE} = 8.2$(s): thus our RB online evaluation is just $2.04\%$ of the FEM computational cost.

\subsection{The Von {K{\'a}rm{\'a}n} plate problem}

We consider now a different problem that can be derived from the classical elasticity equations \cite{ciarlet1988three,ciarlet1997mathematical}. It turns out to be nonlinear and brings with it a lot of technical difficulties. Let us consider an elastic, bidimensional and rectangular plate $\Omega = [0,l] \times [0,1]$ in its undeformed state, subjected to a $\mu$-parametrized external load acting on its edge, then the \emph{Airy stress potential} and the deformation from its flat state, respectively $\phi$ and $u$ are defined by the Von {K{\'a}rm{\'a}n} equations

\begin{equation}
\label{karm}
\begin{cases}
\Delta^2u + \mu u_{xx} = \left[\phi, u\right] + f \ , \quad &\text{in}\ \Omega \\
\Delta^2\phi = -\left[u, u\right] \ , \quad &\text{in}\ \Omega 
\end{cases}
\end{equation}
where $$\Delta^2 := \Delta\Delta = \left(\frac{\p\,^2}{\p\,x^2} + \frac{\p\,^2}{\p\,y^2}\right)^2 \ ,$$ is the biharmonic operator and 
$$[u,\phi] := \frac{\p\,^2u}{\p\,x^2}\frac{\p\,^2\phi}{\p\,y^2} -2\frac{\p\,^2u}{\p x\p y}\frac{\p\,^2\phi}{\p x\p y} + \frac{\p\,^2u}{\p\,y^2}\frac{\p\,^2\phi}{\p\,x^2} \ ,$$
is the \emph{bracket of Monge-Amp\'ere}. So we have a system of two nonlinear and parametrized equations of the fourth order with $\mu$ the parameter that measures the compression along the sides of the plate.
\begin{figure} [htbp]   
\centering    
\input{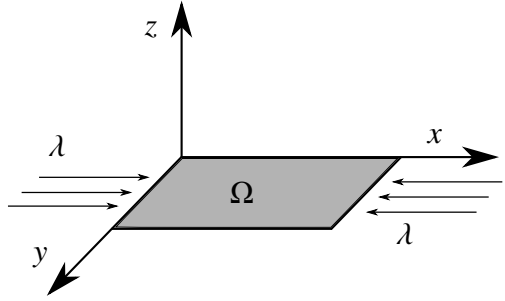}   
\caption{A rectangular bidimensional elastic plate compressed on its edges}  
\label{piastra} 
\end{figure}

From the mathematical point of view, we will suppose the plate is simply supported, i.e. that holds boundary conditions $$u = \Delta u = 0, \qquad \phi = \Delta \phi = 0, \qquad \text{on} \ \p\Omega .$$
In this model problem we are interested in the study of stability and uniqueness of the solution for a given parameter. In fact due to the nonlinearity of the bracket we obtain the so called \emph{buckling phenomena} \cite{PAMM:PAMM201310213}, that is the main feature studied in bifurcations theory. 
What we seek is the critical value of $\mu$ for which the stable (initial configuration) solution become unstable while there are two new stable and symmetric solutions.

To detect this value we need a very complex algorithm that mixes a \emph{continuation method}, a nonlinear solver and finally a full-order method to find the buckled state. At the end for every $\mu \in \calD_{train}$ (a fine discretization of the parameter domain $\calD$) we have a loop due to the nonlinearity, for which at each iteration we have to solve the Finite Element method applied to the weak formulation of the problem.

Here we consider $P=1$ parameter $\mu$ and its domain is suitably chosen\footnote{It is possible to show that the bifurcation point is related to the eigenvalue of the linearized model \cite{berger}, so we are able to set in a proper way the range of the parameter domain.} as $\calD = [30, 70]$.

Also in this case we can simply recover the affine forms with $Q^a = 3$. For the rectangular plate test case with $l=2$ we applied the Finite Element method, with $n_{node} = 441$ nodes and $n_{elem} = 800$ $P_2$ elements, which corresponds to $\calN = 6724$ degrees of freedom. We stress on the fact that the linear system obtained by the Galerkin projection has to be solved at each step of the nonlinear solver, here we chose the classic \emph{Newton method} \cite{QuateroniValli97}.

 Moreover, for a given parameter, we have to solve a FE system until Newton method converges just to obtain one of the possible solutions of our model; keeping in mind that we do not know a priori where is the bifurcation point and we have to investigate the whole parameter domain. It is clear that despite the simple geometry and the quite coarse mesh, the reduction strategies are fundamental in this kind of applications.
 
\begin{figure} [htbp]
\centering
\includegraphics[scale=0.6]{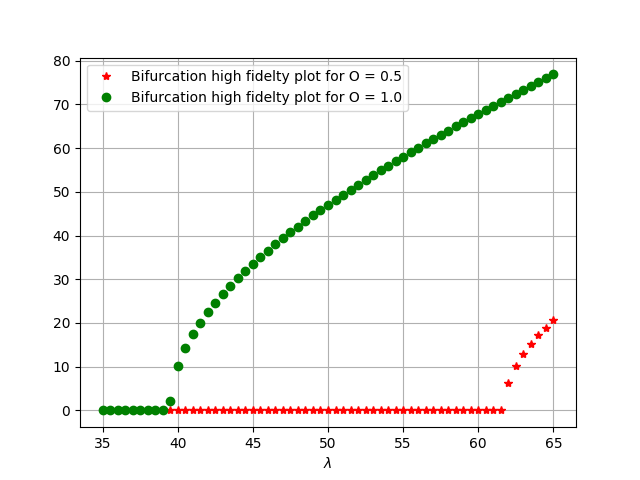}
\caption{Bifurcation diagram for a square plate and different initial guess for Newton method, on y-axis is represented the infinite norm of the solution}
\label{fig:bifdia}
\end{figure}

For example, in order to plot a \emph{bifurcation diagram} like the one in Figure~\ref{fig:bifdia}, the full order code running on  a standart computer takes approximately one hour.

\begin{figure} [htbp]
\centering
\includegraphics[scale=0.4]{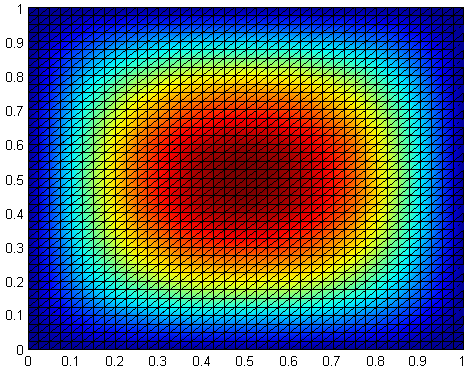} \quad \includegraphics[scale=0.4]{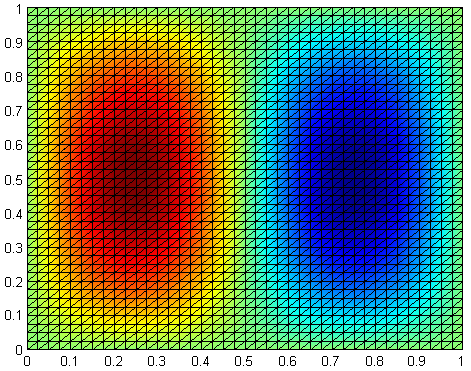}
\caption{Contour plot of the two solutions belonging to the green and red branches of the bifurcation diagram for $\lambda = 70$, respectively}
\label{fig:bifurc1}
\end{figure}

Once selected a specific parameter, $\lambda = 70$, we can see in Figure~\ref{fig:bifurc1} the two solutions that belong to the different branches of the plot reported in Figure~\ref{fig:bifdia}.

 We then applied RB approximation and present in Table~\ref{tab:ex7_tab} a  convergence results: the error between the truth approximation and the reduced one as a function of $N$. The error reported, $\calE_N = \max_{\bfmu \in \calD} ||\bfu^\calN(\bfmu) - \bfu^\calN_{{\rm RB},N}(\bfmu)||_{X}$ is the maximum of the approximation error over a uniformly chosen test sample.

\begin{table}
\centering
\begin{tabular}{|c||c|}
    \hline
		
    $N$ & $\calE_N$  \\
    \hline
    1 &  6.61\texttt{E}+00 \\
    2 &  6.90\texttt{E}-01 \\
 		3 &  7.81\texttt{E}-02 \\
		4 &  2.53\texttt{E}-02 \\
    5 &  1.88\texttt{E}-02 \\
		6 &  1.24\texttt{E}-02 \\
    7 &  9.02\texttt{E}-03 \\
		8 &  8.46\texttt{E}-03 \\
		
    \hline
\end{tabular}
\caption{The Von {K{\'a}rm{\'a}n} plate problem : RB convergence}
\label{tab:ex7_tab}
\end{table}

\begin{figure} [htbp]
\centering
\includegraphics[scale=0.4]{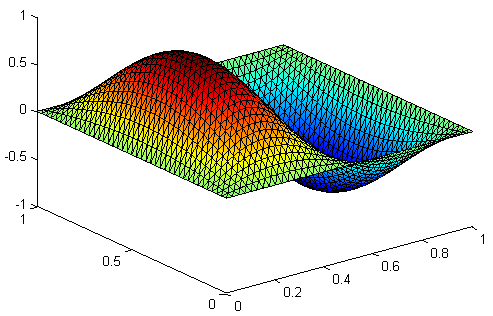} \quad \includegraphics[scale=0.4]{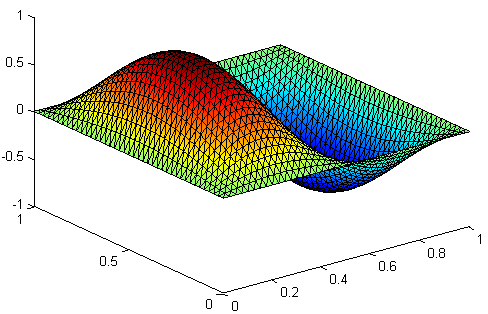}
\caption{Comparison between the full order solution (left) and reduced order one (right) for $\lambda = 65$ }
\label{fig:1cell1}
\end{figure}

As we can see in Figure~\ref{fig:1cell1} e obtain very good results with a low number of snapshots due to the strong properties of the underlying biharmonic operator.

 A suitable extension for the a posteriori error estimate of the solution can be obtained by applying Brezzi-Rappaz-Raviart (BRR) theory on the numerical approximation of nonlinear problems \cite{Brezzi1980,Brezzi1981,Brezzi1982,Grepl2007,canuto2009}. However, the adaptation of BRR theory to RB methods in bifurcating problems is not straightforward, and we leave it for further future investigation \cite{pichirozza}. 

As regards computational times, a RB online evaluation $\bfmu \rightarrow \bfu^\calN_{{\rm RB},N}(\bfmu)$ requires just $t_{\rm RB} = 100$(ms) for $N = 8$; while the FE solution $\bfmu \rightarrow \bfu^\calN(\bfmu)$ requires $t_{\rm FE} = 8.17$(s): thus our RB online evaluation is just $1.22\%$ of the FEM computational cost.

\section{Conclusions}
We have provided some examples of applications of reduced basis
methods in linear elasticity problems depending also on many parameters
of different kind (geometrical, physical, engineering) using
different linear elasticity approximations, a 2D Cartesian setting or a 3D axisymmetric one, different material models (isotropic and orthotropic), as well as an overview on nonlinear problems. Reduced basis methods have confirmed a
very good computational performance with respect to a classical
finite element formulation, not very suitable to solve
parametrized problems in the real-time and many-query contexts.
We have extended and generalized previous work \cite{milani08:RB_LE}  with the possibility to treat with more complex outputs by  introducing a
dual problem \cite{rozza08:ARCME}. Another
very important aspect addressed in this work is the certification of the errors in the
reduced basis approximation by means of a posteriori error
estimators, see for example \cite{huynh07:cras}. This work looks also at more complex 3D parametrized applications (not only in the special axisymmetric case) as quite promising problem to be solved with the same certified methodology \cite{chinesta2014separated,zanon:_phd_thesis}.

\section*{Acknowledgement}
We are sincerely grateful to Prof. A.T. Patera (MIT) and  Dr. C.N. Nguyen (MIT)  for important
suggestions, remarks, insights, and codevelopers of the \texttt{rbMIT}  and \texttt{RBniCS} (http://mathlab.sissa.it/rbnics) software libraries used for the numerical tests. We acknowledge the European Research Council consolidator grant H2020 ERC CoG 2015 AROMA-CFD GA 681447 (PI Prof. G. Rozza).

\section*{Appendix}
\addcontentsline{toc}{section}{Appendix}
\section{Stress-strain matrices}
In this section, we denote $E_i$, $i = 1,3$ as the Young's moduli, $\nu_{ij}$; $i,j = 1,2,3$ as the Poisson ratios; and $G_{12}$ as the shear modulus of the material.
\subsection{Isotropic cases}

For both of the following cases, $E = E_1 = E_2$, and $\nu = \nu_{12} = \nu_{21}$.

Isotropic plane stress:
\begin{equation*}
[\bfE] = \frac{E}{(1-\nu^2)}\left[
         \begin{array}{ccc}
           1 & \nu & 0 \\
           \nu & 1 & 0 \\
           0 & 0 & 2(1+\nu) \\
         \end{array}
       \right].
\end{equation*}

Isotropic plane strain:
\begin{equation*}
[\bfE] = \frac{E}{(1+\nu)(1-2\nu)}\left[
         \begin{array}{ccc}
           1 & \nu & 0 \\
           \nu & 1 & 0 \\
           0 & 0 & 2(1+\nu) \\
         \end{array}
       \right].
\end{equation*}

\subsection{Orthotropic cases}

Here we assume that the orthotropic material axes are aligned with the axes used for the analysis of the structure. If the structural axes are not aligned with the orthotropic material axes, orthotropic material rotation must be rotated by with respect to the structural axes. Assuming the angle between the orthogonal material axes and the structural axes is $\theta$, the stress-strain matrix is given by $[\bfE] = [\bfT(\theta)][\hat{\bfE}][\bfT(\theta)]^T$, where
\begin{equation*}
    [\bfT(\theta)] = \left[
                      \begin{array}{ccc}
                        \cos^2\theta & \sin^2\theta & -2\sin\theta\cos\theta \\
                        \sin^2\theta & \cos^2\theta & 2\sin\theta\cos\theta \\
                        \sin\theta\cos\theta & -sin\theta\cos\theta & \cos^2\theta-\sin^2\theta\\
                      \end{array}
                    \right].
\end{equation*}

Orthotropic plane stress:
\begin{equation*}
[\hat{\bfE}] = \frac{1}{(1-\nu_{12}\nu_{21})}\left[
         \begin{array}{ccc}
           E_1 & \nu_{12}E_1 & 0 \\
           \nu_{21}E_2 & E_2 & 0 \\
           0 & 0 & (1-\nu_{12}\nu_{21})G_{12} \\
         \end{array}
       \right].
\end{equation*}
Note here that the condition
\begin{equation}\label{eqn:ortho_plane_stress}
\nu_{12}E_1 = \nu_{21}E_2
\end{equation}
must be required in order to yield a symmetric $[\bfE]$.

Orthotropic plane strain:
\begin{equation*}
[\hat{\bfE}] = \frac{1}{\Lambda}\left[
         \begin{array}{ccc}
           (1-\nu_{23}\nu_{32})E_1 & (\nu_{12}+\nu_{13}\nu_{32})E_1 & 0 \\
           (\nu_{21}+\nu_{23}\nu_{31})E_2 &(1-\nu_{13}\nu_{31})E_2 & 0 \\
           0 & 0 & \Lambda G_{12} \\
         \end{array}
       \right].
\end{equation*}
Here $\Lambda = (1-\nu_{13}\nu_{31})(1-\nu_{23}\nu_{32})-(\nu_{12}+\nu_{13}\nu_{32})(\nu_{21}+\nu_{23}\nu_{31})$. Furthermore, the following conditions,
$$\nu_{12}E_1 = \nu_{21}E_2, \quad \nu_{13}E_1 = \nu_{31}E_3, \quad \nu_{23}E_2 = \nu_{32}E_3,$$
must be satisfied, which leads to a symmetric $[\bfE]$.

An reasonable good approximation for the shear modulus $G_{12}$ in orthotropic case is given by \cite{carroll98:fem} as
\begin{equation}\label{eqn:ortho_shear_modulus}
    \frac{1}{G_{12}} \approx \frac{(1+\nu_{21})}{E_1} + \frac{(1+\nu_{12})}{E_2}.
\end{equation}

\bibliography{reflib}

\newcommand{\noopsort}[1]{} \newcommand{\printfirst}[2]{#1}
  \newcommand{\singleletter}[1]{#1} \newcommand{\switchargs}[2]{#2#1}
\begin{thebibliography}{10}

\bibitem{almroth78:_autom}
B.~O. Almroth, P.~Stern, and F.~A. Brogan.
\newblock Automatic choice of global shape functions in structural analysis.
\newblock {\em AIAA Journal}, 16:525--528, May 1978.

\bibitem{benner2017model}
P.~Benner, A.~Cohen, M.~Ohlberger, and K.~Willcox.
\newblock {\em Model Reduction and Approximation: Theory and Algorithms}.
\newblock Computational Science and Engineering Series. SIAM, Society for
  Industrial and Applied Mathematics, 2017.

\bibitem{benner2015}
P.~Benner, S.~Gugercin, and K.~Willcox.
\newblock A survey of projection-based model reduction methods for parametric
  dynamical systems.
\newblock {\em SIAM Review}, 57(4):483--531, 2015.

\bibitem{morepas2017}
P.~Benner, M.~Ohlberger, A.~Patera, G.~Rozza, and K.~Urban, editors.
\newblock {\em Model Reduction of Parametrized Systems}.
\newblock Springer International Publishing, 2017.

\bibitem{berger}
M.~Berger.
\newblock On {Von K{\'a}rm{\'a}n} equations and the buckling of a thin elastic
  plate, {I} the clamped plate.
\newblock {\em Communications on Pure and Applied Mathematics}, 20, 1967.

\bibitem{Brezzi1980}
F.~Brezzi, J.~Rappaz, and P.~A. Raviart.
\newblock Finite dimensional approximation of nonlinear problems.
\newblock {\em Numerische Mathematik}, 36(1):1--25, 1980.

\bibitem{Brezzi1981}
F.~Brezzi, J.~Rappaz, and P.~A. Raviart.
\newblock Finite dimensional approximation of nonlinear problems.
\newblock {\em Numerische Mathematik}, 37(1):1--28, 1981.

\bibitem{Brezzi1982}
F.~Brezzi, J.~Rappaz, and P.~A. Raviart.
\newblock Finite dimensional approximation of nonlinear problems.
\newblock {\em Numerische Mathematik}, 38(1):1--30, 1982.

\bibitem{canuto2009}
C.~Canuto, T.~Tonn, and K.~Urban.
\newblock A posteriori error analysis of the reduced basis method for nonaffine
  parametrized nonlinear pdes.
\newblock {\em SIAM Journal on Numerical Analysis}, 47(3):2001--2022, 2009.

\bibitem{carroll98:fem}
W.~F. Carroll.
\newblock {\em A Primer for Finite Elements in Elastic Structures}.
\newblock Wiley, 1998.

\bibitem{MOR2016}
F.~Chinesta, A.~Huerta, G.~Rozza, and K.~Willcox.
\newblock {\em Model Order Reduction: a survey}.
\newblock Wiley Encyclopedia of Computational Mechanics, 2016.

\bibitem{chinesta2014separated}
F.~Chinesta and P.~Ladev{\`e}ze.
\newblock {\em Separated Representations and PGD-Based Model Reduction:
  Fundamentals and Applications}.
\newblock CISM International Centre for Mechanical Sciences. Springer Vienna,
  2014.

\bibitem{ciarlet1988three}
P.~G. Ciarlet.
\newblock {\em Mathematical Elasticity, Volume I: Three-Dimensional
  Elasticity}.
\newblock Elsevier Science, 1988.

\bibitem{ciarlet1997mathematical}
P.~G. Ciarlet.
\newblock {\em Mathematical Elasticity: Volume II: Theory of Plates}.
\newblock Studies in Mathematics and its Applications. Elsevier Science, 1997.

\bibitem{grepl04:_reduc_basis_approx_time_depen}
M.~A. Grepl and A.~T. Patera.
\newblock {\it A Posteriori} error bounds for reduced-basis approximations of
  parametrized parabolic partial differential equations.
\newblock {\em M2AN (Math. Model. Numer. Anal.)}, 39(1):157--181, 2005.

\bibitem{Grepl2007}
M.A. Grepl, Y.~Maday, N.C. Nguyen, and A.T. Patera.
\newblock Efficient reduced-basis treatment of nonaffine and nonlinear partial
  differential equations.
\newblock {\em ESAIM: Mathematical Modelling and Numerical Analysis},
  41(3):575--605, 8 2007.

\bibitem{hesthaven2015certified}
J.S. Hesthaven, G.~Rozza, and B.~Stamm.
\newblock {\em Certified Reduced Basis Methods for Parametrized Partial
  Differential Equations}.
\newblock SpringerBriefs in Mathematics. Springer International Publishing,
  2015.

\bibitem{hutchingson79:fracture}
J.W. Hutchingson.
\newblock {\em Nonlinear Fracture Mechanics}.
\newblock Monograph, Department of Solid Mechanics, Technical University,
  Denmark, 1979.

\bibitem{huynh07:_phd_thesis}
D.~B.~P. Huynh.
\newblock {\em Reduced-Basis Approximation and Applications in Fracture
  Mechanics}.
\newblock PhD thesis, Singapore-MIT Alliance, National University of Singapore,
  2007.

\bibitem{rbMIT_URL}
D.~B.~P. Huynh, N.~C. Nguyen, G.~Rozza, and A.~T. Patera.
\newblock {\em rbMIT Software:
  http://augustine.mit.edu/methodology/methodology\_rbMIT\_System.htm}.
\newblock Copyright MIT, Technology Licensing Office, case 12600, Cambridge,
  MA, 2007-2009.

\bibitem{huynh07:ijnme}
D.~B.~P. Huynh and A.~T. Patera.
\newblock Reduced basis approximation and a posteriori error estimation for
  stress intensity factors.
\newblock {\em International Journal For Numerical Methods In Engineering},
  72(10):1219–--1259, 2007.

\bibitem{huynh07:cras}
D.~B.~P. Huynh, G.~Rozza, S.~Sen, and A.~T. Patera.
\newblock A successive constraint linear optimization method for lower bounds
  of parametric coercivity and inf-sup stability constants.
\newblock {\em CR Acad Sci Paris Series I}, 345:473--478, 2007.

\bibitem{huynh08:infsupLB}
D.B.P. Huynh, D.~J. Knezevic, Y.~Chen, Jan~S. Hesthaven, and A.~T. Patera.
\newblock A natural-norm {S}uccessive {C}onstraint {M}ethod for inf-sup lower
  bounds.
\newblock {\em Computer {M}ethods in {A}pplied {M}echanics and {E}ngineering},
  199(29-32):1963--1975, 2010.

\bibitem{milani08:RB_LE}
R.~Milani, A.~Quarteroni, and G.~Rozza.
\newblock Reduced basis method for linear elasticity problems with many
  parameters.
\newblock {\em Computer Methods in Applied Mechanics and Engineering},
  197:4812--4829, 2008.

\bibitem{murakami01:SIFhandbook}
Y.~Murakami.
\newblock {\em Stress Intensity Factors Handbook}.
\newblock Elsevier, 2001.

\bibitem{NgocCuong2005}
N.C. Nguyen, K.~Veroy, and A.T. Patera.
\newblock {\em Certified Real-Time Solution of Parametrized Partial
  Differential Equations}, pages 1529--1564.
\newblock Springer Netherlands, Dordrecht, 2005.

\bibitem{noor81:_recen}
A.~K. Noor.
\newblock Recent advances in reduction methods for nonlinear problems.
\newblock {\em Comput. Struct.}, 13:31--44, 1981.

\bibitem{noor82}
A.~K. Noor.
\newblock On making large nonlinear problems small.
\newblock {\em Comp. Meth. Appl. Mech. Engrg.}, 34:955--985, 1982.

\bibitem{noor80:_reduc}
A.~K. Noor and J.~M. Peters.
\newblock Reduced basis technique for nonlinear analysis of structures.
\newblock {\em AIAA Journal}, 18(4):455--462, April 1980.

\bibitem{parks77:a_stiff_sif}
D.~M. Parks.
\newblock A stiffness derivative finite element technique for determination of
  crack tip stress intensity factors.
\newblock {\em International Journal of Fracture}, 10(4):487--502, 1974.

\bibitem{patera07:book}
A.T. Patera and G.~Rozza.
\newblock {\em Reduced basis approximation and A posteriori error estimation
  for Parametrized Partial Differential Equation}.
\newblock MIT Pappalardo Monographs in Mechanical Engineering, Copyright MIT
  (2007-2010).
\newblock http://augustine.mit.edu.

\bibitem{pichirozza}
F.~Pichi and G.~Rozza.
\newblock Reduced basis approaches for parametrized bifurcation problems held
  by nonlinear {V}on {K{\'a}rm{\'a}n} equations.
\newblock {\em Submitted}, 2017.

\bibitem{quarteroni2015reduced}
A.~Quarteroni, A.~Manzoni, and F.~Negri.
\newblock {\em Reduced Basis Methods for Partial Differential Equations: An
  Introduction}.
\newblock UNITEXT. Springer International Publishing, 2015.

\bibitem{Quarteroni2011}
A.~Quarteroni, G.~Rozza, and A.~Manzoni.
\newblock Certified reduced basis approximation for parametrized partial
  differential equations and applications.
\newblock {\em Journal of Mathematics in Industry}, 1(1):3, 2011.

\bibitem{QuateroniValli97}
A.~Quateroni and A.~Valli.
\newblock {\em Numerical Approximation of Partial Differential Equations}.
\newblock Springer, 2nd edn., 1997.

\bibitem{rozza08:ARCME}
G.~Rozza, D.B.P. Huynh, and A.T. Patera.
\newblock Reduced basis approximation and a posteriori error estimation for
  affinely parametrized elliptic coercive partial differential equations:
  Application to transport and continuum mechanics.
\newblock {\em Archives of Computational Methods in Engineering}, 15:229--275,
  2008.

\bibitem{CMCS-CONF-2009-002}
G.~Rozza, N.~C. Nguyen, A.T. Patera, and S.~Deparis.
\newblock Reduced {B}asis {M}ethods and {A} {P}osteriori {E}rror {E}stimators
  for {H}eat {T}ransfer {P}roblems.
\newblock In {\em Proceedings of the {ASME} {HT} 2009 summer conference}, Heat
  Transfer Confererence, New York, 2009. ASME.
\newblock Presented in the Computational Section of the conference, Paper No.
  HT2009-88211, pp. 753-762.

\bibitem{sneddon1999mathematical}
I.N. Sneddon, R.~Dautray, P.~Benilan, J.L. Lions, M.~Cessenat, A.~Gervat,
  A.~Kavenoky, and H.~Lanchon.
\newblock {\em Mathematical Analysis and Numerical Methods for Science and
  Technology: Volume 1 Physical Origins and Classical Methods}.
\newblock Mathematical Analysis and Numerical Methods for Science and
  Technology. Springer Berlin Heidelberg, 1999.

\bibitem{sneddon2000mathematical}
I.N. Sneddon, R.~Dautray, P.~Benilan, J.L. Lions, M.~Cessenat, A.~Gervat,
  A.~Kavenoky, and H.~Lanchon.
\newblock {\em Mathematical Analysis and Numerical Methods for Science and
  Technology: Volume 2 Functional and Variational Methods}.
\newblock Mathematical Analysis and Numerical Methods for Science and
  Technology. Springer Berlin Heidelberg, 2000.

\bibitem{veroy03:_phd_thesis}
K.~Veroy.
\newblock {\em Reduced-Basis Methods Applied to Problems in Elasticity:
  Analysis and Applications}.
\newblock PhD thesis, Massachusetts Institute of Technology, 2003.

\bibitem{roark01:roark_formula}
W.~Young and R.~Budynas.
\newblock {\em Roark's Formulas for Stress and Strain}.
\newblock McGraw-Hill Professional, Seventh Edition, 2001.

\bibitem{zanon:_phd_thesis}
L.~Zanon.
\newblock {\em Model Order Reduction for Nonlinear Elasticity: Applications of
  the Reduced Basis Method to Geometrical Nonlinearity and Finite Deformation}.
\newblock PhD thesis, RWTH Aachen University, 2017.

\bibitem{PAMM:PAMM201310213}
L.~Zanon and K.~Veroy-Grepl.
\newblock The reduced basis method for an elastic buckling problem.
\newblock {\em PAMM}, 13(1):439--440, 2013.

\bibitem{zienkiewics05:FEM1}
O.~C. Zienkiewicz, R.~L. Taylor, and J.~Z. Zhu.
\newblock {\em The Finite Element Method: Its Basis and Fundamentals}.
\newblock Butterworth-Heinemann, Sixth Edition, 2005.

\end{thebibliography}
\bibliographystyle{plain}

\end{document}